\documentclass[11pt]{article}

\usepackage{latexsym}
\usepackage{amssymb}
\usepackage{amsthm}
\usepackage{amscd}
\usepackage{amsmath}
\usepackage{mathrsfs}
\usepackage[all]{xy}
\usepackage{graphicx}
\usepackage{hyperref} 
\usepackage{tabmac}

\input xy \xyoption{frame}

\usepackage[usenames,dvipsnames]{color}

\theoremstyle{definition}
\newtheorem* {schensted*}{Schensted's Theorem}

\theoremstyle{definition}
\newtheorem* {theorem*}{Theorem}

\newtheorem{theorem}{Theorem}[section]
\newtheorem{question}[theorem]{Question}
\theoremstyle{definition}

\newtheorem* {example*}{Example}

\newtheorem{lemma}[theorem]{Lemma}
\theoremstyle{definition}
\newtheorem{definition}[theorem]{Definition}
\theoremstyle{definition}

\newtheorem{conjecture}[theorem]{Conjecture}
\newtheorem{proposition}[theorem]{Proposition}
\newtheorem{corollary}[theorem]{Corollary}
\newtheorem* {remark}{Remark}
\theoremstyle{definition}
\newtheorem {example}[theorem]{Example}
\theoremstyle{definition}

\theoremstyle{definition}

\theoremstyle{definition}

\xyoption{dvips}

\usepackage{fullpage}

\numberwithin{equation}{section}

\def\ben{\begin{enumerate}}
\def\een{\end{enumerate}}

\newcommand{\larc}[1]{\hspace{-.4ex}\overset{#1}{\frown}\hspace{-.4ex}}
\def\({\left(}
\def\){\right)}

    \newcommand{\sP}{\Lambda}        \newcommand{\RR}{\mathbb{R}} \newcommand{\QQ}{\mathbb{Q}}   \newcommand{\cP}{\mathcal{P}}

\newcommand{\PP}{\mathbb{P}}

\def\NN{\mathbb{N}}
    \def\ZZ{\mathbb{Z}}  \def\X{\mathfrak{X}}

\newcommand{\hs}{\hspace{0.1mm}}

\def\qquand{\qquad\text{and}\qquad}

\def\sP{\mathscr{P}}

\newcommand{\ba}{\begin{aligned}}
\newcommand{\ea}{\end{aligned}}
\newcommand{\barr}{\begin{array}}
\newcommand{\earr}{\end{array}}
\newcommand{\be}{\begin{equation}}
\newcommand{\ee}{\end{equation}}

\def\cV{\mathcal V}

\makeatletter
\renewcommand{\@makefnmark}{\mbox{\textsuperscript{}}}
\makeatother

\allowdisplaybreaks[1]

\UseCrayolaColors

\begin{document}
\title{Crossings and nestings in colored set partitions}
\author{
Eric Marberg
\\ Department of Mathematics \\ Massachusetts Institute of Technology, United States \\ \tt{emarberg@math.mit.edu}
  }
 \date{}

\maketitle

\def\KK{\mathbb{K}}
\def\omdef{:}
\def\bflambda{\mbox{\boldmath$\lambda$}}
\def\bfmu{\mbox{\boldmath$\mu$}}
\def\bfvarnothing{\mbox{\boldmath$\varnothing$}}
\def\sbflambda{\mbox{\scriptsize\boldmath$\lambda$}}
\def\bfP{\mbox{\boldmath$P$}}
\def\bfQ{\mbox{\boldmath$Q$}}
\def\content{\mathrm{content}\hs}
\def\shape{\mathrm{shape}\hs}

\def\Cr{\mathrm{cr}}
\def\Ne{\mathrm{ne}}

\def\eCr{{\mathrm{cr}}}
\def\eNe{{\mathrm{ne}}}

\def\NCN{\mathrm{NCN}}
\def\NC{\mathrm{NC}}
\def\Arc{\mathrm{Arc}}
\def\YY{\mathbb{Y}}
\def\LL{\mathbb{L}}

\def\cVZ{\mathcal{V}^{2n}_{(\varnothing,\dots,\varnothing)}}
\def\cOZ{\mathcal{O}^{2n}_{(\varnothing,\dots,\varnothing)}}
\def\sP{\Pi}

\def\seq{\mathrm{seq}}
\def\cG{\mathcal{G}}

\def\rsk{\overset{\mathrm{RSK}} \longmapsto}
\newcommand\stirl[2]{\genfrac{\{}{\}}{0pt}{}{#1}{#2}}
\newcommand\C[2]{C_{#1}(#2)}
\newcommand{\cGG}[3]{\cG_{#1,#2,#3}}

\begin{abstract}
Chen, Deng, Du, Stanley, and Yan introduced the notion of $k$-crossings and $k$-nestings for set partitions, and proved that the sizes of the largest $k$-crossings and $k$-nestings in the  partitions of an $n$-set possess a symmetric joint distribution. This work considers a generalization of these results to  set partitions whose arcs are labeled by an $r$-element set (which we call \emph{$r$-colored set partitions}). In this context, a $k$-crossing or $k$-nesting  is   a sequence of arcs, all with the same color,  which form a $k$-crossing or $k$-nesting in the usual sense. After showing that the sizes of the largest crossings and nestings in colored set partitions likewise have a symmetric joint distribution, we consider several related enumeration problems.  We prove that  $r$-colored set partitions with no crossing arcs of the same color are in bijection with certain  paths in  $\NN^r$, generalizing the correspondence between noncrossing (uncolored) set partitions and 2-Motzkin paths. Combining this  with recent work of Bousquet-M\'elou and Mishna affords a proof that the sequence counting noncrossing 2-colored  set partitions is P-recursive. We also  discuss how our methods extend to several variations of colored set partitions with analogous notions of crossings and nestings. \end{abstract}


\section{Introduction and statement of results}

\subsection{Introduction}

A \emph{partition} of a set $S$ is a set of disjoint nonempty sets (here called \emph{blocks}) whose union is $S$. 
Given a  partition $P$ of the set $[n]\omdef=\{1,2,\dots,n\}$, we write $\Arc(P)$ for the set  of pairs of integers $(i,j)$  which occur in the same block of $P$ such that $j$ is the least element of the block greater than $i$.    One often  depicts set partitions  by drawing the graph whose vertex set is $[n]$ and whose edge set is $\Arc(P)$; e.g.,
\[P= \{ \{ 1,3,4,7\} ,\{2,6\}, \{5\} \} \quad\text{is represented by}\quad \xy<0.0cm,-0.0cm> \xymatrix@R=-0.0cm@C=.3cm{
*{\bullet} \ar @/^.7pc/ @{-} [rr]   & 
*{\bullet} \ar @/^1.0pc/ @{-} [rrrr] &
*{\bullet} \ar @/^.5pc/ @{-} [r]  &
*{\bullet} \ar @/^0.7pc/ @{-} [rrr] &
*{\bullet} &
*{\bullet}  &
*{\bullet}\\
1   & 
2 &
3  &
4 &
5 &
6 & 
7
}\endxy.\]
Call this graph  the \emph{standard representation} of $P$.

Such pictures motivate the following terminology, introduced in \cite{Stan}.     A \emph{$k$-crossing} or \emph{$k$-nesting} of a  set partition $P$  is a sequence of arcs $\{(i_t,j_t)\}_{t \in [k]}\subseteq \Arc(P)$ satisfying respectively
\be\label{cross-nest}
\ i_1<i_2<\dots<i_k < j_1 <j_2<\dots<j_k
\qquad\text{or}\qquad i_1<i_2<\dots<i_k < j_k  <\dots<j_2< j_1.
\ee 
In the standard representation of a set partition, $k$-crossings and $k$-nestings appear as follows:
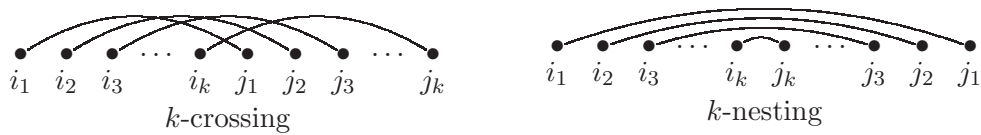
\[
\barr{c} \xy<0.0cm,-0.2cm> \xymatrix@R=-0.0cm@C=0.1cm{
*{\bullet} \ar @/^1.2pc/ @{-} [rrrrr]   & 
*{\bullet} \ar @/^1.2pc/ @{-} [rrrrr] &
*{\bullet} \ar @/^1.2pc/ @{-} [rrrrr] &
*{\dots}  &
*{\bullet} \ar @/^1.2pc/ @{-} [rrrrr] &
*{\bullet} &
*{\bullet}  &
*{\bullet}  &
*{\dots} &
*{\bullet} \\
i_1 & i_2 & i_3 & & i_k & j_1 & j_2 & j_3&& j_k
}\endxy
\\
\text{$k$-crossing}
\earr
\qquad
\barr{c} \xy<0.0cm,-0.0cm> \xymatrix@R=-0.0cm@C=0.1cm{
*{\bullet} \ar @/^1.2pc/ @{-} [rrrrrrrrr]   & 
*{\bullet} \ar @/^.9pc/ @{-} [rrrrrrr] &
*{\bullet} \ar @/^0.6pc/ @{-} [rrrrr] &
*{\dots}  &
*{\bullet} \ar @/^.3pc/ @{-} [r] &
*{\bullet} &
*{\dots}  &
*{\bullet}  &
*{\bullet} &
*{\bullet} \\
i_1 & i_2 & i_3 & & i_k & j_k & & j_3& j_2& j_1
}\endxy
\\
\text{$k$-nesting}
\earr
\]
A partition $P$ of $[n]$ with no $k$-crossing (respectively, $k$-nesting) is \emph{$k$-noncrossing} (respectively, \emph{$k$-nonnesting}), and following \cite{Stan} we let $\Cr(P)$ and $\Ne(P)$ denote the largest integers $k$ such that $P$ has a $k$-crossing or $k$-nesting, respectively.

It is well-known that the  noncrossing (i.e., 2-noncrossing) partitions of $[n]$ and the nonnesting (i.e., 2-nonnesting) partitions of $[n]$ are both counted by the Catalan numbers $C_n \omdef=\frac{1}{n+1} \binom{2n}{n}$.
In \cite{Stan}, Chen {et al.} generalized this fact considerably, proving the following; here $\min(P)$ and $\max(P)$ denote the sets of minimum and maximum elements of the blocks of a set partition $P$.
 \begin{theorem}\label{intro-thm}
 Fix a positive integer $n$ and subsets $S,T \subseteq [n]$.
The statistics $\Cr(P)$ and $\Ne(P)$ have a symmetric joint distribution over all partitions $P$ of $[n]$ with $\min(P) = S$ and $\max(P) = T$.
\end{theorem}

  In other words,
the number of  partitions of $[n]$ which are $j$-noncrossing and $k$-nonnesting is equal to the number of those which are $k$-noncrossing and $j$-nonnesting \cite[Theorem 1.1]{Stan}.
The results of \cite{Stan} have been reinterpreted and extended in a number of ways; we mention without hope of being comprehensive the papers \cite{Burill,BMP,Tangle,DeMier,Jel1,Jel2,Jonsson,Krattenthaler,Rubey,Rubey2,Yen}. Most recently, Chen and Guo \cite{ChenGuo} have generalized the equidistribution of crossings and nestings to colored complete matchings. This paper begins with some enumerative problems which one encounters on extending Chen and Guo's findings   to all {colored set partitions}.



If $r$ is a fixed positive integer, then an \emph{$r$-colored  partition} of $[n]$ is a pair $\Lambda = (P,\varphi)$, consisting of a partition
$P$  of $[n]$ together with a map  $\varphi : \Arc(P) \to [r]$ labeling its arcs.
This fairly natural set partition analogue has appeared most prominently in recent years in the study of the representation theory of  the group 
 of unipotent upper triangular $n\times n$ matrices over a finite field; see \cite{Thiem} for a concise overview of this connection.
%
On the other hand, this notion of a colored set partition dates at least back to \cite{Rogers}, where  
it is studied by Rogers under the name of a ``colored rhyming scheme.''
We also mention that the polynomials counting the $r$-colored partitions of $[n]$ (a variant of the Touchard polynomials) define several  sequences noted in \cite{BerSlo}. 


Following \cite{ChenGuo}, we say that an $r$-colored set partition $\Lambda = (P,\varphi)$
 has a $k$-crossing (respectively, $k$-nesting) if $P$ has a $k$-crossing (respectively, $k$-nesting) involving arcs \emph{all of the same color} with respect to $\varphi$. Define 
$\Cr(\Lambda)$ and $\Ne(\Lambda)$
   as the maximum integers $k$ such that $\Lambda$ has a $k$-crossing or $k$-nesting, respectively. As in the uncolored case we say that  $\Lambda$ of $[n]$ is \emph{$k$-noncrossing} if $\Cr(\Lambda) < k$ and \emph{$k$-nonnesting} if $\Ne(\Lambda) < k$.

It follows as a straightforward corollary of the results in \cite{Stan}  that the joint distribution of the numbers $\Cr(\Lambda)$ and $\Ne(\Lambda)$ over $r$-colored set partitions is also symmetric.
We present the  derivation of this result here in the form of two short lemmas.
First, we note that an $r$-colored set partition  may be viewed as an $r$-tuple of uncolored set partitions satisfying a certain condition:

\begin{lemma}\label{intro-lem}
Given an $r$-colored partition $\Lambda$ of $[n]$,  let 
$\Lambda_t$ for each $ t \in [r]$ be the uncolored partition of $[n]$ for which $\Arc(\Lambda_t)$ is the set of $t$-colored arcs of $\Lambda$. The map 
$\Lambda \mapsto (\Lambda_1,\Lambda_2,\dots, \Lambda_r)$
is then a bijection from $r$-colored partitions of $[n]$ to $r$-tuples $(P_1,\dots,P_r)$ of  uncolored partitions  of $[n]$ with the property that 
$ \min(P_i) \cup \min(P_j) = \max(P_i) \cup \max(P_j) = [n]$ for any distinct $i,j \in [r]$.
\end{lemma}

\begin{proof}
It is easy to see that $r$-colored set partitions are in bijection with $r$-tuples of set partitions whose arc sets have pairwise disjoint  left/right  endpoints. The lemma follows as the sets of right and left endpoints of the arcs of a set partition $P$ are respectively $[n] \setminus\min(P)$ and $[n]\setminus\max(P)$.
\end{proof}

Next, we state a general corollary of Theorem \ref{intro-thm}.
Here we let $\cP([n])$ denote the set of all subsets of $[n]$ and write $\Pi_n$ for the set of uncolored partitions of $[n]$.

\begin{lemma}\label{general-lem} Suppose  $\X$ is a set with an injective map $\phi : \X \to (\Pi_{n})^r$. Let $f$ be any function   with domain $\NN^r$ and define for $x \in \X$ with $\phi(x) = (P_1,\dots,P_r) \in (\Pi_n)^r$
\[  \Cr_{\phi,f}(x) \omdef= f(\Cr(P_1),\dots,\Cr(P_r)) \qquand \Ne_{\phi,f}(x) \omdef = f(\Ne(P_1),\dots,\Ne(P_r)).\]
If  the image of $\X$ under $\phi$ is equal to the inverse image  in $(\Pi_{n})^r$ of some subset of  $ (\cP([n])\times \cP([n]))^{r}$  under the map 
$(P_t)_{t\in [r]} \mapsto \( \min(P_t),\max(P_t)\)_{t \in [r]}$,
then the statistics
$ \Cr_{\phi,f}(x) $ and $ \Ne_{\phi,f}(x)$ for $x \in \X$ 
 possess a symmetric joint distribution.
\end{lemma}

\begin{proof}
Theorem \ref{intro-thm} shows that there exists an involution of $(\Pi_n)^r$ interchanging  the crossing and nesting numbers of an $r$-tuple $(P_1,\dots,P_r)$. The lemma follows since our condition on $\phi$ ensures that $\X$ may be identified with a subset of $(\Pi_n)^r$ which is invariant under this involution.
\end{proof}

The point of these lemmas is the following extension of  \cite[Theorem 1.1]{Stan}. Here, given integers $j,k$ and subsets $S,T \subseteq [n]$,  
 we write $\NCN^{S,T}_{j,k}(n,r)$ for the number of $r$-colored partitions $\Lambda = (P,\varphi)$ of $[n]$ with
$\Cr(\Lambda) < j$ and $\Ne(\Lambda) < k$ and $\min(\Lambda) = S$ and $\max(\Lambda) = T$,
where we define $\min(\Lambda) \omdef = \min(P)$ and $\max(\Lambda)\omdef= \max(P)$.

\begin{theorem}\label{thm-intro}
$\NCN^{S,T}_{j,k}(n,r) = \NCN^{S,T}_{k,j}(n,r)$ for all integers $j,k$ and subsets $S,T \subseteq [n]$. 
\end{theorem}

\begin{proof}
In the notation of Lemma \ref{intro-lem}, we have $\min(\Lambda) = \bigcap_{t\in [r]} \min(\Lambda_t)$ and $\max(\Lambda) = \bigcap_{t \in [r]} \max(\Lambda_t)$. Hence,  the theorem follows  by
applying Lemma \ref{general-lem} with $\phi$ the map in Lemma \ref{intro-lem} and $f(x_1,\dots,x_r) \omdef = \max\{ x_1,\dots,x_r\}$. In particular, the image of $\phi$ on $r$-colored partitions $\Lambda$ with $\min(\Lambda) = S$ and $\max(\Lambda)=T$ is completely characterized by conditions involving only $\min(\cdot)$ and $\max(\cdot)$.
\end{proof}

Let $\NCN_{j,k}(n,r)$ denote the number of all $r$-colored $j$-noncrossing $k$-nonnesting partitions of $[n]$. Summing the previous result over all $S,T\subseteq [n]$ gives the following generalization of \cite[Corollaries 1.2 and 1.3]{Stan}. Here we also write $\NC_k(n,r)$ and $\mathrm{NN}_k(n,r)$ for the number of $r$-colored partitions of $[n]$ which are respectively $k$-noncrossing and $k$-nonnesting.

\begin{corollary}\label{cor-intro} $\NCN_{j,k}(n,r) = \NCN_{k,j}(n,r)$ and $\NC_k(n,r) = \mathrm{NN}_k(n,r)$ for all   $j,k,n,r$.
\end{corollary}


We are thus left with this motivating question: what are the numbers which appear on either side of the equalities in this corollary?
There is an established industry (see, e.g., \cite{BousquetMelou2,JinNC,MansourSeverini,MishnaNN})
 dedicated to producing formulas and generating functions for $\NCN_{j,k}(n,1)$ and $\NC_k(n,1)$ and their variants. From among a plenitude of interesting facts, we mention  that $\{\NCN_{2,2}(n,1)\}_{n=0}^\infty$ gives the sequence odd-indexed Fibonacci numbers \cite[A001519]{OEIS} while $\{\NC_2(n,1)\}_{n=0}^\infty$ gives the sequence of Catalan numbers \cite[A000108]{OEIS}. 
Less is known about these numbers for  values of $r \neq 1$, and this work represents an attempt to begin filling this gap in our understanding.

\subsection{Results}

Our main results  appear in  Sections  \ref{corr-sect}, \ref{enum-sect}, and \ref{extend-sect}, and are summarized as follows.
 In Section \ref{corr-sect} we adapt the results of \cite{Stan} to describe a correspondence between $j$-noncrossing $k$-nonnesting $r$-colored set partitions and walks on a certain multigraph (see Theorem \ref{graph-thm}). Using this 
we are able to show  that for each fixed $j,k,r$, 
the ordinary generating function $\sum_{n\geq 0} \NCN_{j,k}(n,r) x^n $ is a rational  power series (see Corollary \ref{rat-cor}).
%
Moreover, define 
\[\C{n}{r} \omdef= \NC_2(n,r) = \mathrm{NN}_2(n,r)\] as the number of $r$-colored  noncrossing (or nonnesting) partitions of $[n]$,
so that $\C{n}{1}  = \frac{1}{n+1} \binom{2n}{n}$.
After Theorem \ref{graph-thm} below, we  prove the following result, which one can view as a generalization of the well-known fact that the Catalan number $C_n$ is equal to the number of \emph{2-Motzkin paths} of length $n-1$. 
Note here that $\NN^r$ denotes the set of vectors in $\RR^r$ with nonnegative integer coordinates.

\begin{theorem}\label{nc-cor}
Fix  positive integers $n$ and $r$ and let $e_i$ denote the $i{\mathrm{th}}$ unit coordinate vector in $\RR^r$.
Then 
$\C{n}{r}$
 is 
equal to the number of $(n-1)$-step walks in  $\NN^r$ from the origin to itself using the $(r+1)^2$ steps given by 
$\pm e_i$ for $i \in [r]$ (contributing $2r$ steps), $e_i-e_j$ for $i\neq j$ (contributing $r^2-r$ steps), and $r+1$ distinct zero steps.
\end{theorem}



Using this theorem, we can prove a more explicit result concerning the enumeration of noncrossing 2-colored set partitions. The sequence of values of $\C{n}{2}$ begins as follows (and appears now as \cite[A216947]{OEIS}):
\[ \{ \C{n}{2}\}_{n=0}^\infty = \( 1,\ 1,\ 3,\ 11,\ 47,\ 225,\ 1173,\ 6529,\ 38265,\ 233795, \dots\).\]
While these numbers seem not to have a simple closed formula, we can at least establish the following statement. 
There is an interesting parallel between this result 
and \cite[Proposition 1]{BousquetMelou2}, which asserts something  similar for the number of 3-noncrossing (uncolored) partitions of $[n]$. 

\begin{theorem}\label{intro-prop} If $n$ is a positive integer then $\C{n}{2}$ is equal to the constant coefficient of 
\[ \( 1 - x^2 y^{-1} +x^3  - x^2y^2 + y^3 - x^{-1} y^2\) \( 3+x+y+x^{-1} + y^{-1} + xy^{-1} + x^{-1} y\)^{n-1},\]
and the following polynomial recurrence holds for all nonnegative integers $n$:
\[ 9n(n+3)\cdot \C{n}{2}  - 2\bigl(5n^2+26n+30\bigr) \cdot \C{n+1}{2} + (n+4)(n+5) \cdot \C{n+2}{2} = 0.\]
Thus the generating function $\sum_{n\geq 0} \C{n}{2} x^n $ is  D-finite. 
\end{theorem}

Our proof of this result appears at the end of Section \ref{enum-sect},  and combines work of Bousquet-M\'elou and Mishna on walks in the quarter plane \cite{BM} with  Zeilberger's algorithm for creative telescoping \cite[Chapter 6]{Zeil}.

\begin{remark} Standard techniques \cite{Wimp} for determining the asymptotics of solutions to linear recurrence equations with polynomial coefficients show that 
\[\C{n}{2} \sim \kappa\cdot 9^n/n^4\quad\text{as $n\to \infty$}\] for some positive real constant $\kappa$. (We used Zeilberger's {\sc{Maple}} package {\sc{AsyRec}} \cite{asyrec} to derive this growth rate automatically.) 
Empirical estimates (using the {\tt{AsyC}} command in \cite{asyrec}) indicate that
$ \kappa = \frac{3^5}{2^4} \frac{\sqrt{3}}{\pi}$, an equality one can establish rigorously using arguments similar to those  in \cite[\S2.6]{BousquetMelou2}. 
Interestingly, these results show precisely that 
$ C_n(2) \sim 3v_n$ as $n\to \infty$, where $v_n$ is the number of  1234-avoiding permutations of $[n]$ \cite[A005802]{OEIS}.
\end{remark}

Of course the generating function $\sum_{n\geq 0} \C{n}{1} x^n$ is also D-finite, since the Catalan numbers satisfy the polynomial recurrence $2(2n+1)C_n = (n+2) C_{n+1}$ for all $n \in \NN$. One naturally asks whether  $\sum_{n\geq 0} \C{n}{r} x^n$ is likewise D-finite for any integers $r >2$. Empirical evidence suggests a negative answer to this question, at least for $r=3$. One can efficiently compute values of $\C{n}{r}$ using Proposition \ref{nc-cor} with a standard lattice path counting algorithm (which stores intermediate data in an intelligent way to avoid repetitious calculations).
In this way we have computed $\C{n}{3}$ for $n\leq 200$, which gives enough values to detect   a polynomial recurrence of  order $\leq 12$  with coefficients of degree $\leq 12$; however, no such recurrence exists. 

This  mirrors the situation for $k$-noncrossing (uncolored) partitions of $[n]$. In \cite{BousquetMelou2}, Bousquet-M\'elou and  Xin  prove 
that the generating function $\sum_{n\geq 0} \NC_3(n,1) x^n$  is D-finite, but  
identify several reasons  (see \cite[\S4]{BousquetMelou2}) why it seems improbable that $\sum_{n\geq 0} \NC_k(n,1) x^n$ is D-finite for any integers $k\geq 4$.
Following their lead,  
we make this conjecture:
 
\begin{conjecture}
The generating function $\sum_{n\geq 0} \C{n}{r} x^n $ is not D-finite  for integers $r\geq 3$.
\end{conjecture}

Several variations of the results in \cite{Stan} have appeared for other objects besides set partitions for which natural concepts of crossings and nestings exists (see, e.g., \cite{BMP,ChenGuo,Tangle,DeMier,Jonsson,Krattenthaler,Yen}). 
In  Section \ref{extend-sect}, we briefly survey  several colored variations of these results. Our main finding is that many colored objects$-$such as matchings \cite{ChenGuo}, permutations \cite{BMP,Yen}, and tangled diagrams \cite{Tangle}$-$can be realized as colored set partitions $\Lambda$ satisfying certain conditions on $\min(\Lambda)$ and $\max(\Lambda)$; see Propositions \ref{enhanced-part-match-lem}, \ref{permutation-lem}, and \ref{tangled-lem}. The techniques we used to prove Theorem \ref{thm-intro} can therefore also be used to easily derive the symmetric joint distribution of crossing and 
nesting numbers in these cases. Using this idea we are able to generalize and simplify the proofs of some related results on crossings and nestings in the literature.

\subsection*{Acknowledgements}

I am grateful  to Cyril Banderier, Joel Brewster Lewis, Alejandro Morales, Alexander Postnikov, Steven V Sam,  and Richard P. Stanley for  helpful discussions and suggestions.

\section{Preliminaries}\label{prelim-sect}

Here we briefly recollect some of the main results and notation from \cite{Stan}, to be adapted to colored set partitions in the next section. 
Throughout, we let $\NN$ and $\PP$ denote the sets of nonnegative and positive integers.
For us, a \emph{partition} of an integer $n$ is a weakly decreasing sequence of  positive integers $\lambda=(\lambda_1,\lambda_2,\dots,\lambda_\ell)$ with $|\lambda|\omdef = \sum_{i=1}^\ell \lambda_i =n$.  Define $\lambda_i$ to be $0$ for all $i \in \PP$ exceeding $\ell$; 
the Young diagram of $\lambda$ is then the set $\{ (i,j) \in \PP \times \PP : j \leq \lambda_i\}$, which we represent  as a left-justified array of boxes with $\lambda_i$ boxes in row $i$, as in the following example:
\[\lambda=(4,2,1)\quad\text{has Young diagram}\quad {\tiny\tableau[s]{&&& \\ & \\ \ }}.\]
To ``add a box'' to a partition $\lambda$ means to produce a partition $\mu$ whose Young diagram is obtained by  adding a single box to that of $\lambda$.  
Deletion of boxes is defined similarly. Let $\YY$ denote the set of partitions of nonnegative integers; this set is partially ordered by inclusion of Young diagrams.

In the spirit of \cite{Stan}, we adopt the following terminology:

\begin{definition} \label{def1} A \emph{semi-oscillating tableau} (of length $n$) is  a sequence of partitions
 $(\lambda^0,\lambda^1,\dots,\lambda^n)$ with $\lambda^0 = \lambda^n=\varnothing$ such that  $\lambda^{i}$ is obtained  from $\lambda^{i-1}$ for each $i \in [n]$  by either adding a box, deleting a box, or doing nothing (so that $\lambda^i = \lambda^{i-1}$).  
\end{definition} 

For example, $\( 
      \varnothing,\ 
      {\tableau[p]{\ }},\ 
      {\tableau[p]{ \\ \ }},\ 
            {\tableau[p]{ \\ \ }},\ 
      {\tableau[p]{ &  \\   \ }},\ 
      {\tableau[p]{ \\  \ }},\ 
      {\tableau[p]{ \ }},\ 
            {\tableau[p]{ \ }},\ 
      \varnothing \)$ is a semi-oscillating tableau of length $n=8$.
      
      \begin{remark}
Here we insist on the convention $\lambda^n =\varnothing$; note however that the definitions in  \cite{Stan} do not make this requirement when defining various kinds of analogous tableaux.
\end{remark}

      For us, a \emph{matching} is  a set partition whose blocks each have either one or two elements. 
(Sometimes this term refers to set partitions whose blocks all have size two, which we refer to as \emph{complete matchings}.) Given a matching $M$ of $[n]$, let $a_j$ for $j \in [n]$ be the unique number such that $\{j,a_j\}$ is a block of $M$. (Note it is possible to have $a_j=j$.) Define  $\pi^j_M$ for $0\leq j \leq n$ as the subsequence of $a_1a_2\cdots a_j$ with all letters $\leq j$ removed, and let  $\lambda_M^j$ be the integer partition which is the common shape of the pair of SYT's assigned to $\pi_M^j$ by the RSK correspondence, as defined in \cite[\S7.11]{StanleyEnum2}. 

      \begin{example}\label{ex1} For the matching 
\[
   M = \{\{1,8\}, \{2,5\}, \{3\}, \{4,6\}, \{7\} \}= \xy<0.0cm,-0.0cm> \xymatrix@R=-0.0cm@C=.3cm{
*{\bullet} \ar @/^1.2pc/ @{-} [rrrrrrr]   & 
*{\bullet} \ar @/^0.7pc/ @{-} [rrr] &
*{\bullet} &
*{\bullet} \ar @/^0.6pc/ @{-} [rr]  &
*{\bullet}  &
*{\bullet} &
*{\bullet} &
*{\bullet}  \\
1   & 
2 &
3  &
4 &
5 &
6 &
7 &
8
}\endxy
\]
we have $\{ \pi^j_M\}_{0\leq j \leq 8}= \(\emptyset, 8,85,85,856,86,8,8,\emptyset\)$ and 
$\{ \lambda^j_M\}_{0\leq j \leq 8} = \( 
      \varnothing,\ 
      {\tableau[p]{\ }},\ 
      {\tableau[p]{ \\ \ }},\ 
            {\tableau[p]{ \\ \ }},\ 
      {\tableau[p]{ &  \\   \ }},\ 
      {\tableau[p]{ \\  \ }},\ 
      {\tableau[p]{ \ }},\ 
            {\tableau[p]{ \ }},\ 
      \varnothing \).$
      \end{example}

The following simple statement includes several main results in \cite{Stan} specialized to the case of matchings. We will generalize this theorem to colored set partitions by a sequence of brief lemmas, and so have sketched a proof to make everything done here more self-contained.

\begin{theorem}[See Chen \emph{et al.} \cite{Stan}] \label{match-map}
The map
$M \mapsto (\lambda_M^0, \lambda_M^1,\dots,\lambda_M^n)$
   defines  a bijection from the set of matchings of $[n]$ to the set of  semi-oscillating tableaux of length $n$, such that:
   \begin{itemize}
\item[(a)]  $\Cr(M)$ is the maximum number of columns  occurring in any of the integer partitions $\lambda^{i}_M$.

\item[(b)]  $\Ne(M)$ is the maximum number of rows  occurring in any of the integer partitions $\lambda^{i}_M$.

\item[(c)] $\min(M) $ is the set of $i \in [n]$ with $ \lambda^{i-1}_M  \subseteq \lambda^{i}_M$.

\item [(d)] $\max(M) $ is the set of $i \in [n]$ with $\lambda^{i-1}_M \supseteq \lambda^{i}_M$.
   \end{itemize}
   \end{theorem}
   
 \begin{proof}[Proof Sketch]  
 The given map 
  is a bijection  by arguments similar to (and easier than) the proofs of \cite[Theorems 2.4 and 3.2]{Stan}.
  Properties (c) and (d) are immediate,  while properties (a) and (b) hold because
$M$ has a 
  $k$-crossing (respectively, a $k$-nesting) if and only if some sequence $\pi^j_M$ has an increasing (respectively, a decreasing) subsequence of length $k$, which
occurs, by Schensted's theorem \cite{Schensted},  if and only if the partition $\lambda^j_M$  has $k$ columns (respectively, $k$ rows).
\end{proof}

The preceding theorem   extends  from matchings to arbitrary set partitions by the following.  

\begin{lemma}\label{part-match-lem} 
The   
 map sending a partition  $P$ of $[n]$ to the unique matching $M $ of $[2n]$ with 
 \[\Arc(M) = \left\{ (2i,2j-1) : (i,j) \in \Arc(P)\right\}\] is a bijection 
 from the set of  partitions of $[n]$ to the set of matchings $M$ of $[2n]$ such that $2i \in \min(M)$ and $2i-1 \in \max(M)$ for all $i \in [n]$.
Furthermore,  $\Cr(P) = \Cr(M)$ and $\Ne(P) = \Ne(M)$.
 \end{lemma} 
 
 \begin{proof}[Proof Sketch] The lemma is intuitively clear since our map  is  defined by applying the local rules
  \[
  \xy<0.2cm,0.0cm> \xymatrix@R=0.0cm@C=0.2cm{
 &&&& \\
 &&
*{\bullet} \ar @/_.6pc/ @{-} [ull] \ar @/^.6pc/ @{-} [urr]   & &
}\endxy
\mapsto \hspace{0mm}
  \xy<0.2cm,0.0cm> \xymatrix@R=-0.0cm@C=0.2cm{
 &&&&& \\
 &&
*{\bullet} \ar @/_.6pc/ @{-} [ull] 
& *{\bullet}  \ar @/^.6pc/ @{-} [urr]   & &
}\endxy
\quad
  \xy<0.2cm,0.0cm> \xymatrix@R=-0.0cm@C=0.2cm{
 &&& \\
 &&
*{\bullet} \ar @/_.6pc/ @{-} [ull]
}\endxy
\mapsto  \hspace{0mm}
  \xy<0.2cm,0.0cm> \xymatrix@R=-0.0cm@C=0.2cm{
 &&&&& \\
 &&
*{\bullet} \ar @/_.6pc/ @{-} [ull] 
& *{\bullet}     & &
}\endxy
\quad
  \xy<0.2cm,0.0cm> \xymatrix@R=-0.0cm@C=0.2cm{
&& \\
*{\bullet}  \ar @/^.6pc/ @{-} [urr]   & &
}\endxy
\mapsto  \hspace{0mm}
  \xy<0.2cm,0.0cm> \xymatrix@R=-0.0cm@C=0.2cm{
 &&&& \\
 &
*{\bullet}
& *{\bullet} \ar @/^.6pc/ @{-} [urr]     & &
}\endxy
\quad
  \xy<0.2cm,0.0cm> \xymatrix@R=0.0cm@C=0.2cm{
 && \\
 &
*{\bullet}   &
}\endxy
\mapsto  
  \xy<0.2cm,0.0cm> \xymatrix@R=-0.0cm@C=0.2cm{
 &&& \\
 &
*{\bullet} 
& *{\bullet}   & 
}\endxy
\]
 to the standard representation of a set partition $P$ of $[n]$. The details are left to the reader.
\end{proof}

Composing the maps in the preceding theorem and lemma gives a bijection from set partitions of $[n]$ to semi-oscillating tableaux $(\lambda^0,\lambda^1,\dots,\lambda^{2n})$ which are \emph{vacillating} in the following sense.  

\begin{definition} \label{def2}
A \emph{vacillating tableau}  is a semi-oscillating tableau $(\lambda^0,\lambda^1,\dots,\lambda^n)$ which has $\lambda^{i-1} \subseteq \lambda^{i}$ when $i$ is even and $ \lambda^{i-1} \supseteq \lambda^{i}$ when $i$  is odd.
\end{definition}

A vacillating tableau with length $n=11$ is
$\( 
      \varnothing,\ 
            \varnothing,\ 
      {\tableau[p]{\ }},\ 
            {\tableau[p]{\ }},\ 
      {\tableau[p]{ \\ \ }},\ 
            {\tableau[p]{ \\ \ }},\ 
      {\tableau[p]{ &  \\   \ }},\ 
      {\tableau[p]{ \\  \ }},\ 
            {\tableau[p]{ \\  \ }},\ 
      {\tableau[p]{ \ }},\ 
            {\tableau[p]{ \ }},\ 
      \varnothing \)$. The following theorem from \cite{Stan} is immediate from combining Theorem \ref{match-map} and Lemma \ref{part-match-lem}.

\begin{theorem}[See Chen \emph{et al.} \cite{Stan}]\label{sp-version}
There is a bijection from the set of $j$-noncrossing $k$-nonnesting partitions $P$ of $[n]$ to the set of vacillating tableaux $(\lambda^0,\lambda^1,\dots,\lambda^{2n})$ for which every $\lambda^i$ has fewer than $j$ columns and fewer than $k$ rows. Further, this bijection is such that $i \in \min(P)$ if and only if $\lambda^{2i-2} = \lambda^{2i-1}$ and $i \in \max(P)$ if and only if $\lambda^{2i-1}  =\lambda^{2i}$.
\end{theorem}

Theorem \ref{intro-thm}  follows as a corollary of  this result on noting that the component-wise transpose of  integer partitions defines an involution of the set of vacillating tableaux, which interchanges the maximum numbers of rows and columns and which preserves the indices where $\lambda^{i-1} = \lambda^i$.


\section{Colored set partitions and  $r$-partite tableaux}\label{corr-sect}


We now describe one generalization of the results in the previous section which will prove useful in enumerating certain classes of colored set partitions. (Chen and Guo investigate another generalization in \cite{ChenGuo}; the relationship between our statements and those in Chen and Guo's work will be discussed in Section \ref{matching-sect}.) To begin we note the following definition:


\begin{definition} Let $r$ be a positive integer. An \emph{$r$-partite partition of $n$} is a sequence $\bflambda = (\lambda^1,\lambda^2,\dots,\lambda^r)$ of integer partitions such that $\sum_{i=1}^r |\lambda^i| = n$. 
\end{definition}

\begin{remark}
These sequences  provide a common indexing set for the conjugacy classes and irreducible characters of the wreath product $ (\ZZ/r \ZZ) \wr S_n$ of a cyclic group by a symmetric group; see for example \cite[Section 2]{APR2008}.
We will typically use {boldface} symbols to indicate $r$-partite partitions.  
\end{remark}

The Young diagram of an $r$-partite partition is just the sequence of Young diagrams of its components, and with respect to this convention, the addition and deletion of boxes is defined exactly as for ordinary integer partitions.
We  therefore define \emph{semi-oscillating}  and \emph{vacillating $r$-partite tableaux} exactly as in Definitions \ref{def1} and \ref{def2}, only now as sequences of  $r$-partite partitions instead of  integer partitions. 
Let us define also an \emph{oscillating $r$-partite tableaux} to be a semi-oscillating tableaux $(\bflambda^0,\bflambda^1,\dots,\bflambda^n)$ such that $\bflambda^{i-1} \neq \bflambda^i$ for all $i$.

We now have this statement extending Theorem \ref{sp-version}.

\begin{theorem} \label{summary-thm}
For any positive integers $n$ and $r$, there are bijections
\[ \barr{rcl}
\text{$r$-colored partitions of $[n]$} & \leftrightarrow& \text{vacillating $r$-partite tableaux of length $2n$} \\
\text{$r$-colored matchings of $[n]$} & \leftrightarrow &  \text{semi-oscillating $r$-partite tableaux of length $n$}  \\
\text{$r$-colored complete matchings of $[2n]$} & \leftrightarrow &  \text{oscillating $r$-partite tableaux of length $2n$}.
\earr
\]
With respect to each bijection, if $\Lambda \mapsto (\bflambda^0,\bflambda^1,\dots)$ 
then $\Cr(\Lambda) < j$ and $\Ne(\Lambda) < k$ if and only if the $r$ components of each $\bflambda^i$ all have
 fewer than $j$ columns and fewer than $k$ rows.
\end{theorem}

\begin{proof}
Fix an $r$-colored set partition $\Lambda$ of $[n]$. If we first apply the map in Lemma \ref{intro-lem}  to split $\Lambda$ into $r$ uncolored set partitions $\Lambda_1$, \dots $\Lambda_r$, and then apply to each of these components the map in  Theorem \ref{sp-version}, we obtain a matrix of integer partitions 
$\(  \lambda^{i,t}\)$, where $i=0,1,\dots,2n$ and $t=1,2,\dots,r$, 
such that 
(i) each column $(\lambda^{0,t}, \lambda^{1,t},\dots,\lambda^{2n,t})$ is a vacillating tableaux
and
(ii) in each row at most one  $\lambda^{i,t}$ for $t \in [r]$ differs from its predecessor $\lambda^{i-1,t}$.
Indeed, since $(\Lambda_1$, \dots $\Lambda_r)$ if an arbitrary $r$-tuple of uncolored set partitions such that for any distinct $t,t' \in [r]$ each number $i \in [n]$ belongs to $\min(\Lambda_t)$ or $\min(\Lambda_{t'})$ and also to $\max(\Lambda_t)$ or $\max(\Lambda_{t'})$, it follows from Theorem \ref{sp-version} that
  for each $i \in [2n]$  either $\lambda^{i-1,t} = \lambda^{i,t}$ or $\lambda^{i-1,t'} = \lambda^{i,t'}$. Thus at most one of $\lambda^{i,t}$ or $\lambda^{i,t'}$ can differ from $\lambda^{i-1,t}$ or $\lambda^{i-1,t'}$. 
 
 The matrix $\( \lambda^{i,t}\)$ consequently  represents the same data as the vacillating $r$-partite tableaux  $(\bflambda^0,\bflambda^1,\dots,\bflambda^{2n})$ where $\bflambda^i \omdef= ( \lambda^{i,1},\dots,\lambda^{i,r})$, and the correspondence $\Lambda \mapsto (\bflambda^0,\bflambda^1,\dots,\bflambda^{2n})$ gives the first bijection in the theorem. It is clear from Theorem \ref{sp-version} that $\Cr(\Lambda) < j$ and $\Ne(\Lambda) < k$ if and only if the $r$ components of each $\bflambda^i$ all have
 fewer than $j$ columns and fewer than $k$ rows. Moreover, it follows  that  $i \in \min(\Lambda)$ if and only if $\bflambda^{2i-2} = \bflambda^{2i-1}$ and $i \in \max(\Lambda)$ if and only if $\bflambda^{2i-1}  =\bflambda^{2i}$.
 
This last property implies that $r$-colored matchings of $[n]$ (i.e., set partitions $\Lambda$ with $\min(\Lambda)\cup \max(\Lambda) = [n]$) are in bijection with vacillating $r$-partite tableaux $(\bflambda^0,\bflambda^1,\dots,\bflambda^{2n})$ such that $\bflambda^{2i-2} = \bflambda^{2i-1}$ or $\bflambda^{2i-1}  =\bflambda^{2i}$ for each $i \in [n]$.
We obtain the second bijection in the theorem by noting that such tableaux are in bijection with semi-oscillating $r$-partite tableaux of length $n$
via the map $(\bflambda^0,\bflambda^1,\dots,\bflambda^{2n}) \mapsto (\bflambda^0,\bflambda^2,\bflambda^4,\dots,\bflambda^{2n})$. 

The construction of last bijection follows similarly on noting that $r$-colored complete matchings of $[2n]$ are in bijection with vacillating $r$-partite tableaux $(\bflambda^0,\bflambda^1,\dots,\bflambda^{4n})$ such that exactly one of the equalities $\bflambda^{2i-2} = \bflambda^{2i-1}$ or $\bflambda^{2i-1}  =\bflambda^{2i}$ holds for each $i \in [2n]$.
\end{proof}

Tracing through the details of the preceding discussion affords an explicit description of the bijections in Theorem \ref{summary-thm}. One way of stating this (for colored set partitions) goes as follows. Fix an $r$-colored partition $\Lambda$ of $[n]$, and for each $t \in [r]$ and  $k \in [2n]$, let 
\[ a^t_{k} = \begin{cases} j&\text{if $k=2i$ is even and $(i,j,t) \in \Arc(\Lambda)$ for some $j$},\\ 
0&\text{otherwise}.\end{cases}
\]
Define $\pi_{\Lambda}^{k,t}$  as the subsequence of $a^t_1a^t_2\cdots a^t_k$ with all letters  $\leq  \frac{k+1}{2}$ (and in particular all zeros) removed, so that $\pi_{\Lambda}^{0,t}= \pi_{\Lambda}^{2n,t}=\emptyset$ where $\emptyset$ denotes the empty sequence.
Let
$\bflambda_{\Lambda}^k = \( \lambda_{\Lambda}^{k,1}, \lambda_{\Lambda}^{k,2}, \dots, \lambda_{\Lambda}^{k,r}\)$ be
the $r$-partite partition whose $t^{\mathrm{th}}$ component  is 
the common shape $\lambda_{\Lambda}^{k,t}$ of the pair of SYT's assigned to $\pi_{\Lambda}^{k,t}$ by the RSK correspondence. 
It is a straightforward exercise to check that 
the rule
 \be\label{map-def} \Lambda \mapsto  \(\bflambda_\Lambda^0, \bflambda_\Lambda^1,\dots,\bflambda^{2n}_\Lambda\)\ee
coincides with the composition of the maps in Lemma \ref{intro-lem}, Theorem \ref{match-map}, and Lemma \ref{part-match-lem} and so  has  the properties described 
  in Theorem \ref{summary-thm}.

\def\onebox{{\tableau[p]{\ }}}
\def\twobox{{\tableau[p]{ \\ \ }}}
 
 \begin{example}\label{ex2}
For the 2-colored (2-noncrossing 3-nonnesting) partition  
\[ \Lambda = \left\{\{ 1 \larc{1}4 \larc{2} 5\larc{1}8 \}, \{2\larc{2} 6\larc{1}7\}, \{3\}\right\} =
 \xy<0.0cm,-0.0cm> \xymatrix@R=-0.0cm@C=.3cm{
*{\bullet} \ar @/^.7pc/ @{-}^1 [rrr]   & 
*{\bullet} \ar @/^1.2pc/ @{-}^2 [rrrr] &
*{\bullet} & 
*{\bullet} \ar @/^.3pc/ @{-}^2 [r]  &
*{\bullet} \ar @/^1.2pc/ @{-}^1 [rrr] &
*{\bullet} \ar @/^.3pc/ @{-}^1 [r]  &
*{\bullet}  &
*{\bullet}\\
1   & 
2 &
3  &
4 &
5 &
6 & 
7 &
8
}\endxy
\] 
the sequences $\pi_\Lambda^{k,t}$ and $\bflambda_\Lambda^k = (\lambda_\Lambda^{k,1},\lambda_\Lambda^{k,2})$ are given by

\[ \ba \\[-10pt]
&\begin{array} 
{c | lllllllllllllllll}
\hline
k & 0 & 1 &2 & 3 & 4 & 5 & 6 & 7 & 8 & 9 & 10 & 11 & 12 & 13 & 14 & 15 & 16 \\
\hline
\\[-10pt]
\lfloor\tfrac{k+1}{2}\rfloor & 0 & 1 & 1 & 2 & 2 & 3 & 3 & 4 & 4 & 5 & 5 & 6 & 6 & 7 & 7 & 8 & 8 \\
\\[-10pt]
\hline
\\[-10pt]
\pi_{\Lambda}^{k,1} & \emptyset & \emptyset & 4 & 4 & 4 & 4 & 4 & \emptyset & \emptyset & \emptyset & 8 & 8 &  87 & 8 & 8 & \emptyset & \emptyset \\
\pi_{\Lambda}^{k,2}  & \emptyset & \emptyset & \emptyset & \emptyset & 6 & 6 &6 &  6 & 65 & 6 & 6 & \emptyset &  \emptyset & \emptyset & \emptyset & \emptyset & \emptyset 
\\
\\[-10pt]
\hline
\\[-10pt]
\lambda_{\Lambda}^{k,1} & \varnothing & \varnothing & \onebox & \onebox & \onebox & \onebox & \onebox &  \varnothing & \varnothing & \varnothing & \onebox & \onebox &  \twobox & \onebox & \onebox & \varnothing & \varnothing \\
\lambda_{\Lambda}^{k,2}  & \varnothing & \varnothing & \varnothing & \varnothing & \onebox & \onebox &\onebox &  \onebox & \twobox & \onebox & \onebox & \varnothing &  \varnothing & \varnothing & \varnothing & \varnothing & \varnothing 
\earr
\\
\\[-10pt]
\ea
\]
 \end{example}

We may prove the rationality of some power series by 
translating Theorem \ref{summary-thm} into a statement concerning a bijection between colored set partitions and walks on a certain multigraph (i.e., an undirected graph with multiple edges and loops allowed).
Fix positive integers $j,k,r$.
Since only the trivial  partition of $[n]$ into $n$ blocks is 1-noncrossing or 1-nonnesting, let us assume $j,k\geq 2$.
Now let $\cGG{j}{k}{r}$ denote the multigraph whose 
vertices consist of all $r \times (k-1)$ integer matrices $A = (A_{i\ell}) $ with 
\[ 
j >A_{i,1} \geq  A_{i,2} \geq  \cdots \geq A_{i,k-1} \geq 0
\quad\text{for each }i \in [r],
\]
and which has $e(A,A')$ undirected edges connecting any matrices $A$ and $A'$, where
\begin{itemize}
\item $e(A,A') =  1$ if $A-A'=\pm E_{i\ell}$ or $A-A' = E_{i\ell} - E_{i'\ell'}$ for some  $(i,\ell) \neq (i',\ell')$, where $E_{i\ell}$ denotes the  $r\times (k-1)$ matrix with 1 in position $(i,\ell)$ and 0 in all other positions.

\item $e(A,A') = 1 +d_1+d_2+\dots +d_r$ if $A=A'$, where $d_i$ is the number of distinct entries  in the $i{\mathrm{th}}$ row of $A$ which are less than $j-1$;
\item $e(A,A') = 0$ in all other cases.
\end{itemize}
The following statement  generalizes \cite[Theorem 3.6]{Stan}. 

\begin{theorem}\label{graph-thm} For any positive integers $j,k,r,n$ with $j,k\geq 2$, the number
$\NCN_{j,k}(n,r)$ of $r$-colored $j$-noncrossing $k$-nonnesting partitions of $[n]$  is equal to the number of $(n-1)$-step walks on the multigraph $\cGG{j}{k}{r}$ which begin and end at the zero matrix. \end{theorem}

\begin{proof}
First, observe that we can identify the set of  $r$-partite partitions $\bflambda=(\lambda^1,\lambda^2,\dots,\lambda^r)$  involving only  integer partitions with fewer than $j$ columns and fewer than $k$ rows   with the vertices of $\cGG{j}{k}{r}$ by viewing $\bflambda$ as the $r\times (k-1)$ integer matrix whose $i{\mathrm{th}}$ row list the parts of  $\lambda^i$, possibly extended by zeros.

Now, given a vacillating $r$-partite tableau $  (\bflambda^0,\bflambda^1,\dots,\bflambda^{2n})$ for which each $\bflambda^i$ is a vertex of $\cGG{j}{k}{r}$, consider the sequence of $n$ vertices in $\cGG{j}{k}{r}$ given by $\bflambda^1,\bflambda^3,\bflambda^5,\dots,\bflambda^{2n-1}$. 
Since $n$ is positive this sequence begins and ends at the $r$-partite empty shape, 
 and for each $i \in [n-1]$,  $\bflambda^{2i+1}$ is obtained from $\bflambda^{2i-1}$ by either (1) adding one box, (2) deleting one box, (3) adding one box then deleting one box, or (4) doing nothing.

Cases (1) and (2) correspond in $\cGG{j}{k}{r}$ to adding or subtracting an elementary matrix,
and if $\bflambda^{2i+1} \neq \bflambda^{2i-1}$ then case (3) corresponds in $\cGG{j}{k}{r}$ to adding the difference of two distinct elementary matrices. Thus  if 
$\bflambda^{2i+1} \neq \bflambda^{2i-1}$
then there is exactly one edge between the pair of associated vertices in $\cGG{j}{k}{r}$, and it follows by the definition of a vacillating tableaux that   $\bflambda^{2i}$ is uniquely determined.
On the other hand, the number of ways  
in which one could add a box to $\bflambda^{2i-1}$ without introducing a $j{\mathrm{th}}$ column or $k{\mathrm{th}}$ row and then delete the same box to obtain $\bflambda^{2i+1} = \bflambda^{2i-1}$
is precisely the sum of the numbers of distinct entries less than $k-1$ in the rows of the $r\times (k-1)$ matrix associated to $\bflambda^{2i-1}$. It follows that if $\bflambda^{2i+1} = \bflambda^{2i-1}$ then the number of allowable choices for $\bflambda^{2i}$ is the number of self-loops at the associated vertex in $\cGG{j}{k}{r}$, which suffices to complete the proof of the theorem. 
\end{proof}

\begin{example}\label{graph-ex}
The multigraphs $\cGG{j}{k}{r}$ for $(j,k,r)  \in\{ (2,2,1),(2,2,2)\}$ are shown below, with the vertex corresponding to the zero matrix marked by $\bullet$:
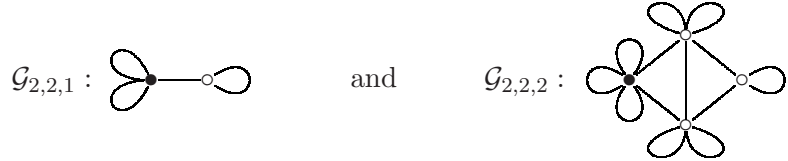
\[
\cGG{2}{2}{1}:
 \xygraph{ 
!{<0cm,0cm>;<1.5cm,0cm>:<0cm,1.2cm>::} 
!~-{@{-}@[|(2.5)]} 
!{(0,0) }*{\bullet}="a" 
!{(0.5,0) }*{\circ}="b" 
"a" -@`{c+(-0.25,-0.75),c+(-0.4,+0.75)} "a" 
"a" -@`{c+(-0.25,0.75),c+(-0.4,-0.75)} "a" 
"b" -@`{p+(0.5,0.5),p+(0.5,-0.5)} "b" 
"a" -@`{p} "b" 
}
\quad\qquad\text{and}\quad\qquad
\cGG{2}{2}{2}:
\xygraph{ 
!{<0cm,0cm>;<1.5cm,0cm>:<0cm,1.2cm>::} 
!~-{@{-}@[|(2.5)]} 
!{(0,0) }*{\bullet}="a" 
!{(0.5,0.5) }*{\circ}="b" 
!{(1.0,0.0) }*{\circ}="c" 
!{(0.5,-0.5) }*{\circ}="d" 
"b" -@`{c+(-0.75,0.25),c+(+0.75,0.4)} "b" 
"b" -@`{c+(0.75,0.25),c+(-0.75,0.4)} "b" 
"d" -@`{c+(-0.75,-0.25),c+(+0.75,-0.4)} "d" 
"d" -@`{c+(0.75,-0.25),c+(-0.75,-0.4)} "d" 
"a" -@ `{p+(0.4,-0.6),p+(-0.4,-0.6)} "a"
"a" -@ `{p+(0.4,0.6),p+(-0.4,0.6)} "a"
"a" -@ `{p+(-0.5,0.5),p+(-0.5,-0.5)} "a"
"c" -@`{p+(0.5,0.5),p+(0.5,-0.5)} "c" 
"a" -@`{p} "b" 
"a" -@`{p} "d"
"c" -@`{p} "b"  
"c" -@`{p} "d" 
"d" -@`{p} "b" 
}
\]
We may explain these pictures by noting that the vertices of $\cG_{2,2,1}$ are the matrices
$A = \(\barr{c} 0 \earr\)$ and $A' = \(\barr{c} 1 \earr\)$
and we have $e(A,A) = 2$ and $e(A,A') = e(A',A') = 1$. Similarly, the vertices of $\cG_{2,2,2}$ are the matrices 
\[B = \(\barr{c} 0 \\ 0 \earr\),\qquad B' = \(\barr{c} 1 \\ 0  \earr\), \qquad B'' = \(\barr{c} 1 \\ 1  \earr\),\qquad\text{and}\qquad B''' = \(\barr{c} 0 \\ 1  \earr\).
\]
One checks that $e(B,B) = 3$ and $e(B',B') = e(B''',B''') = 2$ and $e(B'',B'') = 1$, and that all pairs of distinct matrices are adjacent except $B$ and $B''$.
\end{example}

We now have a short proof of Theorem \ref{nc-cor} from the introduction.

\begin{proof}[Proof of Theorem \ref{nc-cor}]
Fix $n$. For  sufficiently large $j$, we have
 $\C{n}{r} = \NCN_{j,2}(n,r)$ and 
  we may identify 
the   $(n-1)$-step paths in $\cG_{j,2,r}$ from the zero matrix to itself with the $(n-1)$-step paths in $\NN^r$ from the origin to itself 
using the  following $(r+1)^2$ steps:
$\pm e_i$ for $i \in [r]$ (contributing $2r$ steps), or $e_i-e_j$ for $i\neq j$ (contributing $r^2-r$ steps), or $r+1$ distinct zero steps. The number of such paths is therefore equal to $\C{n}{r}$ by Theorem \ref{graph-thm}.
\end{proof}

\section{Enumerating noncrossing colored set partitions}
\label{enum-sect}

We now prove a few enumerative results concerning the numbers $\NCN^{S,T}_{j,k}(n,r)$ (as defined before Theorem \ref{thm-intro}), in particular Theorem \ref{intro-prop} from the introduction. 
To begin, 
 recall the following standard terminology.
Let $\KK$ be a field of characteristic zero.
A formal power series
$y \in \KK[[x]]$  is 
\begin{itemize}
\item \emph{rational} if there are polynomials $P,Q \in \KK[x]$ with $Q(0)\neq 0$ such that $y = PQ^{-1}$ in $\KK[[x]]$;
\item \emph{algebraic} if the powers $1,y,y^2,\dots$ span a finite-dimensional subspace of the set  of formal Laurent series $\KK((x))$, viewed as  a vector space over  the field of rational functions $\KK(x)$;
\item \emph{D-finite} if the derivatives $y,y',y'',\dots$ span a finite-dimensional subspace of $\KK((x))$  over $\KK(x)$.
\end{itemize}
A rational power series is algebraic and an algebraic power series is D-finite, and each class of power series forms a $\KK$-subalgebra of $\KK[[x]]$.
It is useful to note that $y = \sum_{n\geq 0} a_n x^n$ is D-finite if and only if the sequence of coefficients $\{a_n\}$ is \emph{P-recursive}, meaning that 
there exist finitely many polynomials $P_0(x),\dots,P_d(x) \in \KK[x]$, with $P_d(x)$ not identically zero, such that 
\[ P_d(n) a_{n+d} + P_{d-1}(n) a_{n+d-1} + \dots + P_0(n) a_n =0,\qquad\text{for all }n \in \NN.\]  
See \cite[Chapters 4]{StanleyEnum1} and \cite[Chapter 6]{StanleyEnum2} as well as \cite{Lipshitz0,Lipshitz} for more extensive properties of these power series.

Since $\NCN_{j,k}(n,r)$ counts the number of walks of a given length in some graph by Theorem \ref{graph-thm}, the transfer matrix method (see \cite[\S4.7]{StanleyEnum1}) provides an explicit formula for its generating function, and in particular we have this corollary:

\begin{corollary}\label{rat-cor}
For any $j,k,r \in \PP$ the formal power series
 $\sum_{n\geq 0} \NCN_{j,k}(n,r) x^n $ is  rational.
\end{corollary}

\begin{example}

Using \cite[Theorem 4.7.2]{StanleyEnum1} to compute the generating functions for the number of $n$-step walks on the graphs in Example \ref{graph-ex}, we obtain 
\[\sum_{n\geq 0} \NCN_{2,2}(n) x^n = \frac{1-2x}{1-3x+x^2}
\qquad\text{and}\qquad
\sum_{n\geq 0} \NCN_{2,2}(n,2) x^n = \frac{1-6x+7x^2}{1-7x+11x^2-x^3}.\]
Notably, one can derive from this that the number of partitions of $[n]$ which are both noncrossing and nonnesting is $\NCN_{2,2}(n) = f_{2n-1}$ for all $n>0$, where $\{f_n\}_{n=0}^\infty = (0, 1, 1, 2, 3, 5, 8, 13,\dots)$ is the sequence of Fibonacci numbers.  The sequence $\{ \NCN_{2,2}(n)\}_{n=0}^\infty$ appears as \cite[A001519]{OEIS}. More generally, Mansour and Severini have computed an explicit formula for the ordinary generating function of $\{ \NCN_{2,k}(n)\}_{n=0}^\infty$ for any $k\geq 2$ \cite[Theorem 1.1]{MansourSeverini}; it would be interesting to see this result extended to describe the ordinary generating function of  $\{ \NCN_{2,k}(n,r)\}_{n=0}^\infty$.
\end{example}


We might as well also note here another basic property of the numbers $\NCN_{j,k}(n,r)$, which follows  directly from the definitions of $j$-noncrossing and $k$-nonnesting.

\begin{proposition}\label{poly-prop} If $j,k,n$ are fixed positive integers and $S,T \subset [n]$, then   
$\NCN^{S,T}_{j,k}(n,r)$ is a polynomial in $r$  with integer coefficients.
\end{proposition}

\begin{proof}
Fix a partition $P$ of $[n]$ and let $V_{j,k}(P) = \{ \text{$j$-subsets of $\Arc(P)$ which form a $j$-crossing}\} \cup \{ \text{$k$-subsets of $\Arc(P)$ which form a $k$-nesting}\}$.
For each subset $S\subset V_{j,k}(P)$, define   $\cG(S)$ as the graph with  vertex set $S$ which has an edge from $s \in S$ to $s' \in S$ if and only if $s \neq s'$ and $s \cap s' \neq \varnothing$. 
Now let $e(S)$ be the number of connected components of $\cG(S)$ minus $|\bigcup S|$, the number of distinct arcs occurring in elements of $S$.
By the inclusion-exclusion principle, the number of $r$-colorings of $P$ which have no $j$-crossings or  $k$-nestings involving arcs all of the same color is the following sum over all subsets $S$ of $V_{j,k}(P)$:
\[ \sum_{S\subset V_{j,k}(P)} (-1)^{|S|} r^{|\Arc(P)| + e(S)}.\] In particular, this follows since $r^{|\Arc(P)| + e(S)}$ is precisely the number of $r$-colorings of $P$ which assign the same color to every arc in each of the $j$-crossing or $k$-nesting sets $s \in S$.  
Summing this polynomial expression over all $P \in \sP_n$ with $\min(P) = S$ and $\max(P)=T$ gives $\NCN_{j,k}^{S,T}(n,r)$.
\end{proof}

Summing the polynomials $\NCN_{j,k}^{S,T}(n,r)$ over all $S,T\subset[n]$ gives this corollary. 

\begin{corollary}\label{poly-cor} When $j,k\geq 2$, the quantities $\NCN_{j,k}(n,r) = \NCN_{k,j}(n,r)$ and 
$\NC_{k}(n,r) =\mathrm{NN}_k(n,r)$ are monic polynomials in $r$  of degree $n-1$ with integer coefficients.
\end{corollary}

\begin{proof}
The unique partition of $[n]$ with one block and $n-1$ arcs contributes the leading term of $r^{n-1}$ to each of these polynomials when $j,k\geq 2$. 
\end{proof}

The corollary requires $j,k\geq 2$ because if  $k=1$ (or $j=1$) then $\NCN_{j,k}(n,r) = \NC_k(n,r) = 1$ as both polynomials count only the unique partition of $[n]$ with $n$ blocks and no arcs.

\begin{example} 
As in the introduction let  $\C{n}{r} \omdef= \NC_2(n,r)$ be the number of $r$-colored noncrossing partitions of $[n]$. By the corollary this is a polynomial in $r$ of degree $n-1$. The proof of Proposition \ref{poly-prop} outlines an algorithm for computing these polynomials, which we have employed to calculate $\C{n}{r}$ for $n\leq 8$:
    \[
\ba
\C{1}{r}&=  1       \\
\C{2}{r}&=  1 + r   \\      
\C{3}{r}&=  1 + 3r + r^2         \\
\C{4}{r}&=  1 + 5r + 7r^2 + r^3 \\
\C{5}{r}&=  1 + 6r + 19r^2 + 15r^3 + r^4\\
\C{6}{r}&=  1 + 10r + 22r^2 + 67r^3 + 31r^4 + r^5 \\
\C{7}{r}&=  1 + 12r + 56r^2 + 67r^3 + 229r^4 + 63r^5 + r^6 \\
\C{8}{r}&=  1 - 24r + 176r^2 + 159r^3 + 225r^4 + 765r^5 + 127r^6 + r^7.
\ea\]
The pattern displayed by first seven lines of this computation ends at $n=8$, where we see that the coefficients of $\C{n}{r}$ can be  both positive and negative integers.
\end{example}



Before continuing we note  the following  lemma, which will allow us to prove by hand that  the generating function $\sum_{n\geq 0 } \C{n}{2} x^n$ in Theorem \ref{intro-prop} is D-finite. 
This statement is a direct consequence of the main result of \cite{Lipshitz0}, but we have included a proof for completeness.

\begin{lemma} \label{laurent-lem} Assume $\KK$ is a field of characteristic zero and let $S(x) \in \KK\left[x_1^{\pm1},x_2^{\pm1},\dots,x_k^{\pm1}\right]$ be a Laurent polynomial in $k$ indeterminates. 
Fix  $e_1,e_2,\dots,e_k \in \ZZ$
and  define $a_n \in \KK$  for  $n \in \NN$
as the coefficient of $x_1^{e_1}x_2^{e_2}\cdots x_k^{e_k}$ in $S(x)^n$. Then
 the power series $\sum_{n\geq 0} a_n t^n \in \KK[[t]]$
is D-finite.
         \end{lemma}
         
         \begin{remark} If $k=1$ then $\sum_{n\geq 0} a_n t^n$ is in fact algebraic by \cite[Exercise 6.10]{StanleyEnum2}.
         \end{remark}
         
 \begin{proof}
Fix $e=(e_1,e_2,\dots,e_k) \in \ZZ^k$ and define $a_n$ as in the statement of the lemma. Since $a_n$ is unchanged if we simultaneously replace $e_i$ by $-e_i$ and $S(x_1,\dots,x_i,\dots,x_n)$ by $S(x_1,\dots,x_i^{-1},\dots,x_n)$, it is no loss of generality to assume $-e \in \NN^k$. 
 
 Adopt the  notation $x^n \omdef=  x_1^{n_1} x_2^{n_2}\cdots x_k^{n_k}$ for $n=(n_1,n_2,\dots,n_k) \in \ZZ^k$ and identify $\ZZ$ with the diagonal subset of elements $(d,d,\dots,d) \in \ZZ^k$. Also, write $q_n(j)$ for the coefficient of $x^j$ in $S(x)^n$ for $j \in \ZZ^k$, so that $q_n(e) = a_n$.
 
Choose an integer $d \in \NN$  such that $x^dS(x)$ is a polynomial in $\KK[x_1,x_2,\dots,x_k]$, and note that we then have
$q_n(j-dn) = 0$ for all $n \in \NN$ and $ j \in \ZZ^k \setminus \NN^k$. Since $-e \in \NN^k$, so that $e+j \in \NN^k$ only if $j \in \NN^k$,  the following identity of formal power series in $\KK[[x_1,x_2,\dots,x_k,t]]$ thus holds:
 \[ \sum_{j \in \NN^k} \sum_{n \in \NN} q_n(e+ j-dn) x^j t^{dn}
 =
 \sum_{n \in \NN}  x^{dn-e} S(x)^n t^{dn}
 =
    \frac{x^{-e}}{1-x^dt^d S(x) }.
 \] 
In particular, the leftmost power series in $\KK[[x_1,x_2,\dots,x_k,t]]$ is  rational,  and hence D-finite in the sense of \cite{Lipshitz0}. Observe, however, that the diagonal of this power series is precisely $\sum_{n \in \NN} q_n(e) t^{dn} = \sum_{n \in \NN} a_n t^{dn} \in \KK[[t]]$.   \cite[Theorem 1]{Lipshitz0}  asserts that this power series in one variable is   D-finite. This suffices to prove the lemma, finally, because a power series 
 $F(t) \in \KK[[t]]$ is clearly P-recursive (and hence D-finite) if $F(t^d)$ is P-recursive for some positive integer $d$. 
 \end{proof}

We conclude with the proof of Theorem \ref{intro-prop} from the introduction. 
\def\CT{\mathrm{CT}}

\begin{proof}[Proof of Theorem \ref{intro-prop}]
Let $w_n(i,j)$  denote the number of all $n$-step walks in $\ZZ^2$  from $(0,0)$ to $(i,j)$
 using the six steps $\pm(1,0)$, $\pm(0,1)$,  $\pm(1,-1)$,
 and let $q_n(i,j)$ denote the number of such walks which remain in the quarter plane $\NN^2$.
It follows by Theorem \ref{nc-cor} that $\C{n+1}{2} = \sum_{k=0}^n \binom{n}{k} 3^{n-k} q_{k}(0,0)$, and thus the first assertion in the proposition is equivalent to the claim that 
\be\label{fact}q_n(0,0)\quad\text{is the constant term of}\quad\( 1 - x^2 y^{-1} +x^3  - x^2y^2 + y^3 - x^{-1} y^2\)S(x,y)^n,\ee
where $S(x,y)$ is the Laurent polynomial 
$  x+y+x^{-1} + y^{-1} + xy^{-1} + x^{-1} y.$
This fact happens to follow directly from \cite[Proposition 10]{BM}, which asserts more strongly that $q_n(i,j)$ is equal to the coefficient of $x^{-i}y^{-j}$ in the right expression in \eqref{fact} for any $i,j,n \in \NN$.

Bousquet-M\'elou and Mishna prove \cite[Proposition 10]{BM} using the kernel method, a general purpose algebraic argument.  Alternatively, however, one can establish just \eqref{fact} directly from a ``generalized reflection principle'' due to Gessel and Zeilberger \cite{GZ}, in the following way.
Observe that 
\eqref{fact} is  equivalent to the identity
\be\label{gz-id} q_n(0,0)= w_n(0,0) - w_n(-2,1) + w_n(-3,0) - w_n(-2,-2) + w_n(0,-3) - w_n(1,-2).\ee
We will realize this  as a special case of \cite[Theorem 1]{GZ}. Let $R = \{ \pm(e_i-e_j) \in \RR^3 : 1\leq i < j \leq 3\}$ be the root system of type $A_2$  (inheriting the standard inner product $(\cdot,\cdot)$ on $\RR^3$), where $e_i$ denotes the $i$th unit coordinate vector, and fix a  set of simple roots $\Delta  = \{ \alpha,\beta\}$, for example with $\alpha = e_1-e_2$ and $\beta = e_2-e_3$. Define $L$ as the lattice $ \ZZ\alpha\oplus \ZZ\beta$ and let $S$ be the set of six steps 
\[S = \{ \pm (\alpha+\gamma), \pm (\beta+\gamma), \pm(\alpha-\beta)\}\qquad\text{where }\gamma \omdef =\alpha+\beta.\]
The Weyl group $W \cong S_3$ of $R$ preserves both $L$ and $S$. Furthermore, one can show that 
$w_n(i,j)$ is  the number of all $n$-step walks in $L$ from $3\gamma$ to $3\gamma + i(\alpha+\gamma) + j(\beta+\gamma)$  using the steps in $S$, and that $q_n(i,j)$ is the number of such walks 
which stay inside the fundamental Weyl chamber $C \omdef = \{ x \in L : (x,\alpha) > 0\text{ and }(x,\beta) > 0 \}$.
(These statements become clear if one draws  $\alpha $ and $\beta$ as vectors  in the plane they span 
and then works out the subsets of $\RR^2$  corresponding to  $L$, $S$, and $C$; see Figure \ref{fig1}.) 
Now, for these particular choices of $R$, $\Delta$, $L$,  $S$,
the hypotheses of \cite[Theorem 1]{GZ} hold and that theorem (with $a=b=3\gamma$) asserts precisely the identity (\ref{gz-id}).
For a more detailed discussion of this sort of argument, see also 
Grabiner and Magyar's paper \cite{GM}.

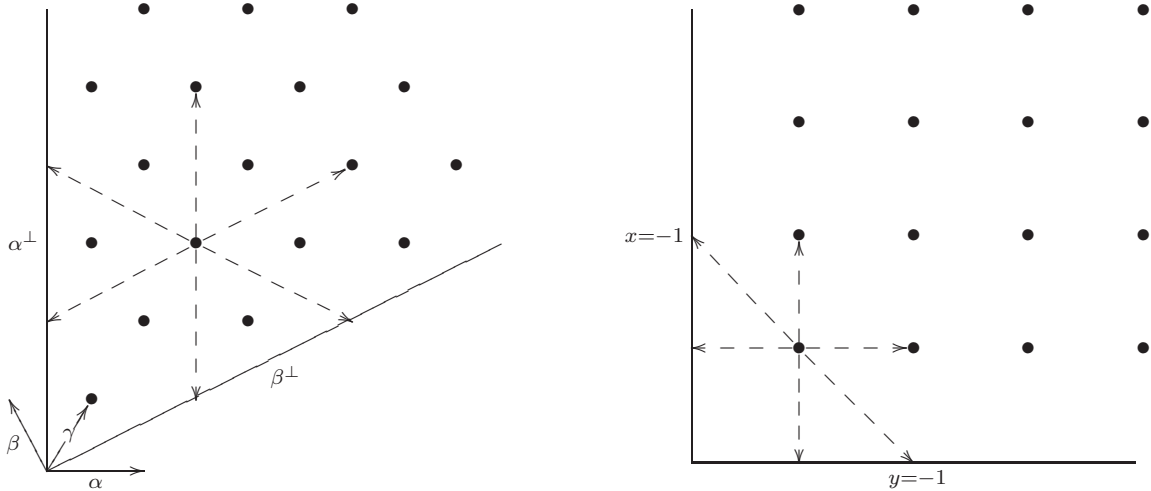
\begin{figure}[h]
\[
\xy<0.0cm,0.0cm> \xymatrix@R=0.8660254037844386cm@C=0.5cm{
*{           }  & *{           } & *{           } & *{\bullet}  & *{           } & *{\bullet} & *{           }  & *{\bullet}  & *{           } & *{           } \\
*{           }  & *{           } & *{\bullet} & *{           }  & *{\bullet} & *{           } & *{\bullet}  & *{           }  & *{\bullet} & *{           } & *{           }\\
*{           }  & *{           } & *{           } & *{\bullet}  & *{           } & *{\bullet} & *{           }  & *{\bullet}  & *{           } & *{\bullet} \\
*{           }  & *{           } & *{\bullet} & *{           }  & 
*{\bullet} 
\ar @{-->} [uu]
\ar @{-->} [dd]
\ar @{-->} [lllu]
\ar @{-->} [rrrd]
\ar @{-->} [rrru]
\ar @{-->} [llld]
& *{           } & *{\bullet}  & *{           }  & *{\bullet} & *{           } & *{           } \\
*{           }  & *{           } & *{           } & *{\bullet}  & *{           } & *{\bullet} & *{           }  & *{           }  & *{           } & *{           } \\
*{           }  & *{           } & *{\bullet} & *{           }  & *{           } & *{           } & *{           }  & *{           }  & *{} & *{           } \\
*{           }  & *{           } \ar @{->} [ul]^\beta  \ar @{->} [rr]_\alpha \ar @{->} [ur]|{\gamma}  
 \ar @{-} [uuuuuu]^{\alpha^\perp}  \ar @{-} [uuurrrrrrrrr]_{\beta^\perp} 
 & *{           } & *{           }  & *{           } & *{} & *{           }  & *{}  & *{           } & *{} \\
}\endxy
\qquad\qquad
 \xy<0.0cm,0.0cm> \xymatrix@R=1.33cm@C=1.33cm{
*{}  & *{\bullet} & *{\bullet} & *{\bullet} & *{\bullet}   \\
*{}  & *{\bullet} & *{\bullet} & *{\bullet}  & *{\bullet}  \\
*{}  & *{\bullet} & *{\bullet} & *{\bullet}  & *{\bullet}  \\
*{}  & *{\bullet} 
\ar @{-->} [u]
\ar @{-->} [d]
\ar @{-->} [r]
\ar @{-->} [l]
\ar @{-->} [ul]
\ar @{-->} [dr]
& *{\bullet} & *{\bullet}  & *{\bullet}  \\
*{}
\ar @{-} [uuuu]^{x=-1}
\ar @{-} [rrrr]_{y=-1}
  & *{} & *{} & *{}   &
}
\endxy
\]
\caption{Comparing walks in $C$ and in $\NN^2$. On the left, the $\bullet$'s are the lattice points in the intersection $C$ of $L$ and the fundamental Weyl chamber of $R$; the solid lines form the boundary of the fundamental chamber; and the dashed arrows, all emanating from the point $3\gamma$,  are the six  steps in $S$. On the right, the $\bullet$'s are the integer lattice points in $\NN^2$; the solid lines are $x=-1$ and $y=-1$; and the dashed lines, all emanating from the origin, are the six steps $\pm(1,0)$, $\pm(0,1)$,  $\pm(1,-1)$.
It is intuitively clear that ``compressing'' and ``dilating'' these pictures defines a bijection between  the walks in $C$ using the steps on the left and the walks in $\NN^2$ using  the steps on the right.
}
\label{fig1}
\end{figure}

It follows from \eqref{fact} and Lemma \ref{laurent-lem} that the ordinary generating function of $\{q_n(0,0)\}_{n=0}^\infty$ is D-finite, so the exponential generating function $Q(x) \omdef= \sum_{n\geq 0} q_n(0,0) \frac{x^n}{n!}$ is D-finite also (as a result of the the equivalence between D-finiteness and P-recursiveness). Let $F(x) = \sum_{n\geq 0}  \C{n}{2} \frac{x^n}{n!}$. The formula $\C{n+1}{2} = \sum_{k=0}^n \binom{n}{k} 3^{n-k} q_{k}(0,0)$ shows that $\frac{d}{dx} F(x) = e^{3x} Q(x)$, so since $e^{3x}$ is D-finite it follows that the derivative $\frac{d}{dx} F(x)$ is D-finite whence $F(x)$ is also D-finite. 
The ordinary generating function $\sum_{n\geq 0} \C{n}{2} x^n$, finally, is therefore  D-finite too. 

To derive the particular polynomial recurrence given in the proposition, we resort to computer methods. Our argument resembles the one employed by Bousquet-M\'elou and Xin to prove \cite[Proposition 1]{BousquetMelou2}. As a preliminary,  observe that 
\[ S(x,y) +2 = 
 y^{-1}(1+x) \Bigl( 1 + x^{-1}(1+x)  y + x^{-1}  y^2\Bigr)
\]
and hence, writing $\CT(f(x,y))$ for the constant term of a polynomial $f(x,y)$, we can directly compute (by considering  first the contribution of the $x^{-1}y^2$ term, then the $x^{-1}(1+x)y$ term, then the remaining terms involving only $x$)
\be\label{ct} \ba  \CT \( x^i y^j (S(x,y)+2)^n\) &= \sum_{k \in \ZZ} \binom{n}{k} \CT\( x^{i-k} y^{j+2k-n} \cdot  (1+x)^n\cdot   \(1 + x^{-1}(1+x)  y\)^{n-k}\)
\\
&
= \sum_{k \in \ZZ} \binom{n}{k} \binom{n-k}{n-j-2k} \CT\( x^{i+j+k-n} \cdot (1+x)^{2n-j-2k}\)
\\
&
=
\sum_{k \in \ZZ} \binom{n}{k} \binom{n-k}{n-j-2k} \binom{2n-j-2k}{n-i-j-k}.
\ea
\ee
%
%
Let $b_n$ denote the constant term of $\( 1 - x^2 y^{-1} +x^3  - x^2y^2 + y^3 - x^{-1} y^2\)(S(x,y)+2)^n$ for $n \in \NN$, and note that $\C{n+1}{2} = \sum_{k=0}^n \binom{n}{k} b_k$.
Equation \eqref{ct} allows us to write $b_n$ as a sum of expressions given by fractions of factorials, to which we can  apply Zeilberger's algorithm for creative telescoping \cite[Chapter 6]{Zeil}. In particular, using the {\sc{Maple}} package {\sc{Ekhad}} \cite{ekhad}, 
we obtain the following 
polynomial recurrence for the numbers $b_n$:
\[ 
(n+5)(n+6)b_{n+2} = 8(n+2)(n+1)b_n+(7n^2+49n+82)b_{n+1},\qquad\text{for all }n \in \NN.
\]
Let $B(x) = \sum_{n\geq 0} b_n \frac{x^n}{n!} \in \QQ[[x]]$ so that
$ e^x B(x) = \sum_{n\geq 0} \C{n+1}{2} \frac{x^n}{n!}$.
It is a routine exercise to check that the preceding recurrence is equivalent to the statement that the formal power series $B(x)$ lies in the kernel of the differential operator
\[ \( \tfrac{d^4}{dx^4} - 7 \tfrac{d^3}{dx^3} - 8\tfrac{d^2}{dx^2}\) x^2 + \( 4\tfrac{d^3}{dx^3}-14 \tfrac{d^2}{dx^2}\) x +\(6\tfrac{d^2}{dx^2}-12\tfrac{d}{dx}\).
\]
Since $(\frac{d}{dx}-1) e^xA(x) = e^x \frac{d}{dx}A(x)$ for any power series $A(x) \in \QQ[[x]]$, 
it follows  that $e^x B(x)$ lies in the kernel of the differential operator
\[ \(\tfrac{d^4}{dx^4} -11 \tfrac{d^3}{dx^3}+19\tfrac{d^2}{dx^2} -9\tfrac{d}{dx}\) x^2 + \( 4\tfrac{d^3}{dx^3}-26\tfrac{d^2}{dx^2}+40\tfrac{d}{dx}-18\) x +\(6 \tfrac{d^2}{dx^2}-24 \tfrac{d}{dx}+18\) ,\]
which in turn implies that the following recurrence holds for all $n \in \NN$:
\[
 \ba 0&= 9n(n+3)\C{n}{2} - (96+97n+19n^2)\C{n+1}{2} +(142+81n+11n^2)\C{n+2}{2} \\&\quad-  (n+5)(n+6)\C{n+3}{2}  . \ea
\]
This recurrence, with the initial conditions $\C{0}{2} = \C{1}{2} = 1$ and $\C{2}{2} = 3$, uniquely determines the sequence $\{ \C{n}{2}\}_{n=0}^\infty$, and it is easy to check that the unique sequence with the same initial conditions  satisfying the three-term polynomial recurrence in the proposition statement  also satisfies this four-term recurrence.
\end{proof}

\section{Some extensions}
\label{extend-sect}

In this final section we discuss a few  variations of  set partitions  with natural notions of crossings, nestings, and colorings. Our discussion here is partially expository, surveying  results from \cite{BMP, ChenGuo,Tangle,Tangle2,Yen}. In each case our noncrossing and nonnesting colored objects  are in bijection with certain classes of noncrossing and nonnesting colored set partitions. On noting these bijections, we are often able to derive the symmetric joint distribution of crossing and nesting numbers by applying the following variant of Lemma \ref{general-lem}.
(Here, as in the introduction, we let $\cP([n])$ denote the set of all subsets of $[n]$. We also write $\Pi_{n,r}$ for the set of $r$-colored partitions of $[n]$.)

\begin{lemma}\label{general-cor} 
Suppose  $\X$ is a set with an injective map $\phi : \X \to \Pi_{n,r}$.
If  the image of $\X$ under $\phi$ is equal to the inverse image  in $\Pi_{n,r}$ of some subset of  $ \cP([n])\times \cP([n])$  under the map 
$\Lambda \mapsto \(\min(\Lambda),\max(\Lambda)\)$,
then the statistics
$\Cr(\phi(x)) $ and $ \Ne(\phi(x))$ for $x \in \X$
 possess a symmetric joint distribution.
\end{lemma}

\begin{proof}
Our argument is almost identical to the short proof of Lemma \ref{general-lem}. Theorem \ref{intro-thm} shows that there exists an involution of $\Pi_{n,r}$ interchanging  all crossing and nesting numbers, and our condition on $\phi$ ensures that $\X$ may be identified with a subset of $\Pi_{n,r}$ which is invariant under this involution. The lemma therefore follows.
\end{proof}


\def\row{\mathrm{row}}
\def\col{\mathrm{col}}

\subsection{Matchings} \label{matching-sect}

A colored set partition $\Lambda$ of $[n]$ is a matching if and only if $\min(\Lambda) \cup \max(\Lambda) = [n]$ (and  a complete matching if and only if furthermore $\min(\Lambda) \cap \max(\Lambda) = \varnothing$). Thus,  
it is immediate from Lemma \ref{general-cor}, taking $\phi$ to be the natural inclusion map, that 
the statistics $\Cr(M)$ and $\Ne(M)$ have a symmetric joint distribution over   $r$-colored  matchings (respectively, complete matchings) $M$ of $[n]$. 
This result is due originally to Chen and Guo \cite{ChenGuo}, who derive it in a considerably different way, by generalizing Theorem \ref{match-map} to the colored matchings using \emph{$r$-rim hook tableaux} in place of $r$-partite tableaux. We briefly review their methods here and explain how we may recover the main results in \cite{ChenGuo} from what is done here by applying a theorem of Fomin and Stanton \cite{FS}.

To this end, recall that if $\mu \subseteq \lambda$ are integer partitions then  the skew shape $\lambda /\mu$ is an \emph{$r$-rim hook} if the skew diagram of $\lambda/\mu$ 
consists of $r$ contiguous squares located on distinct diagonals.  
Let $RH_r$ denote the set of integer partitions $\lambda$ for which there exists a sequence of partitions 
\[\varnothing = \lambda^0 \subset \lambda^1 \subset \dots \subset \lambda^k = \lambda\] such that $\lambda^i / \lambda^{i-1}$ is an $r$-rim hook for each $i \in [k]$.  
The set  $RH_r$, with respect to the partial order in which $\lambda$ covers $\mu$ if and only if  $\mu\subseteq \lambda$ and $\lambda / \mu$ is an $r$-rim hook, is called the \emph{$r$-rim hook lattice}.

Chen and Guo \cite{ChenGuo} define  an \emph{oscillating $r$-rim hook tableau}  as a sequence $(\lambda^0,\lambda^1,\dots,\lambda^n)$ of integer partitions $\lambda^i \in RH_r$ such that $\lambda^0 = \lambda^n = \varnothing$ and either $\lambda^i$ covers $\lambda^{i-1}$ or $\lambda^{i-1}$ covers $\lambda^i$ in $RH_r$ for each $i \in [n]$.
Let us define a \emph{semi-oscillating $r$-rim hook tableau} as a sequence satisfying the same conditions except that it is also allowed that $\lambda^{i-1} = \lambda^i$, and in turn, we define a \emph{vacillating $r$-rim hook tableau} as a semi-oscillating $r$-rim hook tableau  $(\lambda^0,\lambda^1,\dots,\lambda^n)$ with $\lambda^{i-1} \subseteq \lambda^i$ when $i$ is even and $\lambda^{i-1}\supseteq \lambda^i$ when $i$ is odd.
For example,
\[
T = \( 
      \varnothing,\ 
            \varnothing,\ 
      {\tableau[p]{ \\ \ }},\ 
            {\tableau[p]{ \\ \ }},\ 
      {\tableau[p]{ & & \\ \ }},\ 
            {\tableau[p]{ && \\ \ }},\ 
      {\tableau[p]{ & & \\  && }},\ 
      {\tableau[p]{ & \\ & }},\ 
            {\tableau[p]{ & \\  & }},\ 
      {\tableau[p]{ & }},\ 
            {\tableau[p]{ & }},\ 
      \varnothing \)
\]
      is a vacillating 2-rim hook tableau of length $n=11$. 
   These objects are in natural bijection with oscillating, semi-oscillating, and vacillating $r$-partite tableaux by the following  theorem of Fomin and Stanton.
    %
      %
Here, we  view the cartesian product of $r$ copies of the Young lattice $\YY^r$ as the lattice of $r$-partite partitions in which  $\bflambda$  covers $\bfmu$ if and only if $\bflambda$ is obtained from $\bfmu$ by adding one square. 

\begin{theorem}[Fomin and Stanton \cite{FS}] \label{FS-thm}
There is a lattice isomorphism $ RH_r \cong \YY^r$ such that if $\mu \in RH_r$ has $k$ rows and $\ell$ columns and   $\bflambda \in \YY^r$ is the associated $r$-partite partition, then the maximum numbers of rows and columns in the components of $\bflambda$ are 
$\lceil \frac{k}{r}\rceil$  and $\lceil \frac{\ell}{r} \rceil$, respectively.
 \end{theorem}

\begin{proof}
Fomin and Stanton describe a lattice isomorphism $\YY^r \cong RH_r$ in the proof of \cite[Theorem 1.2]{FS}. Given an integer partition $\lambda$, let $f_\lambda : \ZZ \to \ZZ_{\geq 0}$ be the function whose value at $i \in \ZZ$ is the number of squares in the Young diagram of $\lambda$ on the $i{\mathrm{th}}$ diagonal.
If $\bflambda = (\lambda_1,\dots,\lambda_r) \in \YY^r$ and  $\mu \in RH_r$ is the corresponding partition under Fomin and Stanton's isomorphism, then 
$ f_\mu(i) = \sum_{k=1}^r f_{\lambda_k} \( \left\lfloor \tfrac{i + (k-1)}{r} \right\rfloor \)$ for $i \in \ZZ$ \cite[Definition 2.4]{FS}.  Noting this formula, the second assertion in the theorem  follows from that fact that the number of rows (respectively, columns) in a nonempty partition $\lambda$ is the maximum value of $i$ for which $f_\lambda(-i+1)$ (respectively, $f_\lambda(i-1)$) is nonzero.
\end{proof}


Given a semi-oscillating $r$-rim hook tableau $T$, let $\row(T)$ and $\col(T)$ denote the respective maximum number of rows and columns in any partition in $RH_r$ occurring in $T$. 
The following corollary 
 slightly generalizes \cite[Theorem 3.1]{ChenGuo}, which Chen and Guo prove by explicitly describing   a  bijection from $r$-colored complete matchings to oscillating $r$-rim hook tableaux  in terms of 
 Stanton and White's Schensted algorithm for rim hook tableaux \cite{SW}.
 
\begin{corollary}
For any positive integers $n$ and $r$, there are bijections
\[ \barr{rcl}
\text{$r$-colored partitions of $[n]$} & \leftrightarrow& \text{vacillating $r$-rim hook tableaux of length $2n$} \\
\text{$r$-colored matchings of $[n]$} & \leftrightarrow &  \text{semi-oscillating $r$-rim hook tableaux of length $n$}  \\
\text{$r$-colored complete matchings of $[2n]$} & \leftrightarrow &  \text{oscillating $r$-rim hook tableaux of length $2n$}  
\earr
\]
and with respect to each bijection, if $\Lambda \mapsto T$ then $\Cr(\Lambda) = \lceil \col(T) / r \rceil$ and $\Ne(\Lambda) = \lceil \row(T) /r \rceil$.
\end{corollary}

\begin{proof}
Compose the maps in Theorem \ref{summary-thm} with the  bijection between  $r$-partite tableaux and $r$-rim hook tableaux afforded by Theorem \ref{FS-thm}.  
\end{proof}

Let $\NCN^{Matching}_{j,k}(n,r)$ (respectively, $\NC^{Matching}_k(n,r)$) denote the number of $j$-noncrossing $k$-nonnesting (respectively, $k$-noncrossing) $r$-colored tangled diagrams on $[n]$.
The following corollary follows by arguments  similar to the proofs of Theorem \ref{graph-thm} and Corollary \ref{rat-cor}, which we have left to the reader.

\begin{corollary}\label{matching-rat}
For all $j,k,r \in \PP$ the formal power series $\sum_{n\geq 0} \NCN^{Matching}_{j,k}(n,r) x^n$ is rational.
\end{corollary}

In turn we have this corollary:

\begin{corollary}[See Chen and Guo \cite{ChenGuo}]
For all $k,r \in \PP$, $\sum_{n\geq 0} \NC^{Matching}_k(n,r) x^n$ is D-finite.
\end{corollary}

\begin{proof}
Chen and Guo derive an explicit formula \cite[Corollary 3.3]{ChenGuo} for the exponential generating function $ \sum_{n\geq 0} a_n \frac{x^n}{n!}$ where $a_n$ is the number of $r$-colored complete matchings of $[n]$, as the $r$th power of the determinant of a matrix whose entries are linear combinations   of hyperbolic Bessel functions of the first kind $I_n(2x) \omdef = \sum_{i\geq 0} \frac{x^{n+2i}}{i!(n+i)!}$. As $I_n(2x)$ is D-finite and D-finite power series  form a ring,  $\sum_{n\geq 0} a_n \frac{x^n}{n!}$ is D-finite. It follows that  $\sum_{n\geq 0} \NC^{Matching}_k(n,r) \frac{x^n}{n!} = e^x \( \sum_{n\geq 0} a_n \frac{x^n}{n!}\)$ is D-finite, so $\sum_{n\geq 0} \NC^{Matching}_k(n,r) x^n$ is D-finite.
\end{proof}

\subsection{Enhanced set partitions}

In this section we examine a colored version of the enhanced crossing and nesting statistics considered in \cite[Section 4]{Stan}.
Given a partition $P$ of $[n]$, let 
\[\overline{\Arc}(P) = \Arc(P) \cup \{ (i,i) : i \in \min(P) \cap \max(P) \}\] denote the set of \emph{enhanced arcs} of $P$. In turn, define 
an \emph{$r$-colored enhanced partition} of $[n]$ to be a pair $\overline\Lambda = (P,\varphi)$ consisting of a partition $P$ of $[n]$ together with a map $\varphi:\overline{\Arc}(P) \to [r]$. 

The standard representation of an enhanced colored partition is  drawn exactly as for uncolored set partitions, except that we include loops at isolated points (that is, elements of $\min(P) \cap \max(P)$) and label all arcs as in the following example:
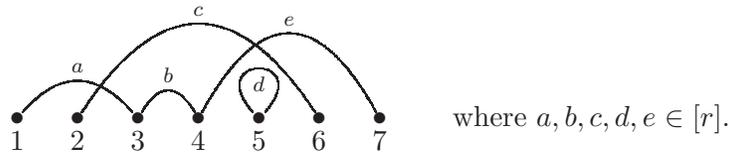
\[ 
 \xy<0.3cm,0.0cm> \xymatrix@R=0.0cm@C=0.4cm{
*{\bullet} \ar @/^1.2pc/ @{-} [rr]^a   & 
*{\bullet} \ar @/^3.0pc/ @{-} [rrrr]^c &
*{\bullet} \ar @/^0.9pc/ @{-} [r]^b  &
*{\bullet} \ar @/^2.7pc/ @{-} [rrr]^e &
*{\bullet} \ar @(ul,ur) @{-}_d &
*{\bullet}  &
*{\bullet}\\
1   & 
2 &
3  &
4 &
5 &
6 & 
7
}\endxy\qquad\text{where }a,b,c,d,e \in [r].
 \]
Let 
$\min(\overline \Lambda) = \min(P)$
and 
$\max(\overline \Lambda) = \max(P)$
and say that $\overline \Lambda$ 
has
an \emph{enhanced $k$-crossing} or  \emph{enhanced $k$-nesting} if there is a sequence    of $k$ enhanced arcs $\{(i_t,j_t)\}_{t \in [k]}\subset \overline{\Arc}(P)$, all labeled with the same color by $\varphi$, satisfying respectively
\be\label{enhanced-cross-nest}
\ i_1<i_2<\dots<i_k \leq j_1 <j_2<\dots<j_k
\qquad\text{or}\qquad i_1<i_2<\dots<i_k \leq j_k  <\dots<j_2< j_1.
\ee
Define  $\eCr(\overline \Lambda)$ and $\eNe(\overline \Lambda)$ as the largest integers $k$ such that $\overline \Lambda$ has an enhanced $k$-crossing or enhanced $k$-nesting, respectively. 
Observe that 1-colored enhanced set partitions represent the same data as uncolored set partitions; in the case $r=1$, the enhanced crossing and nesting statistics given here coincide with those defined for (uncolored, unenhanced) set partitions in \cite{Stan}.

We now have this variant of Lemma \ref{part-match-lem}:

\begin{proposition}\label{enhanced-part-match-lem}
The map which sends  an $r$-colored enhanced partition $\overline \Lambda = (P,\varphi)$ to the unique $r$-colored matching $M = (P',\varphi')$ of $[2n]$ such that 
\begin{itemize}
\item[] $\Arc(P') = \Bigl\{ (2i-1,2j) : (i,j) \in \overline{\Arc}(P)\Bigr\}$ and $\varphi'(2i-1,2j) = \varphi(i,j)$ for $(i,j) \in \overline{\Arc}(P)$
\end{itemize}
is a bijection from  $r$-colored enhanced partitions of $[n]$ 
  to $r$-colored partitions $M$ of $[2n]$ such that $2i-1 \in \min(M)$ and $2i \in \max(M)$ for all $i \in [n]$. Furthermore $\eCr(\overline\Lambda) = \Cr(M)$ and $\eNe(\overline\Lambda) = \Ne(M)$. 
 \end{proposition}
 
 \begin{proof}[Proof Sketch]
 Like Lemma \ref{part-match-lem}, this result is intuitively clear since our map defines $M$ as the colored matching whose standard representation is formed by applying the local rules
  
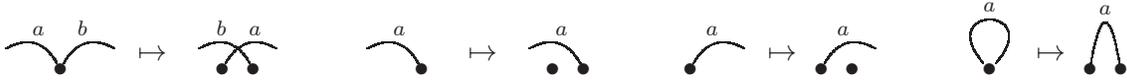
\[
  \xy<0.2cm,0.0cm> \xymatrix@R=0.0cm@C=0.2cm{
 &&&& \\
 &&
*{\bullet} \ar @/_.6pc/ @{-} [ull]_a \ar @/^.6pc/ @{-} [urr]^b   & &
}\endxy
\mapsto \hspace{-3mm}
  \xy<0.2cm,0.0cm> \xymatrix@R=-0.0cm@C=0.2cm{
 &&&&& \\
 &&
*{\bullet} \ar @/^.6pc/ @{-} [urr]^{\hspace{4mm}a}  
& *{\bullet}  \ar @/_.6pc/ @{-} [ull]_{\hspace{-2mm}b}   & &
}\endxy
\quad
  \xy<0.2cm,0.0cm> \xymatrix@R=-0.0cm@C=0.2cm{
 &&& \\
 &&
*{\bullet} \ar @/_.6pc/ @{-} [ull]_a
}\endxy
\mapsto  \hspace{-3mm}
  \xy<0.2cm,0.0cm> \xymatrix@R=-0.0cm@C=0.2cm{
 &&&&& \\
 &&
*{\bullet}
& *{\bullet}  \ar @/_.6pc/ @{-} [ull]_a     & &
}\endxy
\quad
  \xy<0.2cm,0.0cm> \xymatrix@R=-0.0cm@C=0.2cm{
&& \\
*{\bullet}  \ar @/^.6pc/ @{-} [urr]^a   & &
}\endxy
\mapsto  \hspace{-3mm}
  \xy<0.2cm,0.0cm> \xymatrix@R=-0.0cm@C=0.2cm{
 &&&& \\
 &
*{\bullet}\ar @/^.6pc/ @{-} [urr]^a
& *{\bullet}     & &
}\endxy
\quad
  \xy<0.2cm,0.0cm> \xymatrix@R=0.0cm@C=0.2cm{
 && \\
 &
*{\bullet}   \ar@(ul,ur) @{-}^a&
}\endxy
\mapsto  \hspace{-3mm}
  \xy<0.2cm,0.0cm> \xymatrix@R=-0.0cm@C=0.2cm{
 &&& \\
 &
*{\bullet}  \ar @/^1.5pc/ @{-} [r]^a
& *{\bullet}   & 
}\endxy
\]
to the standard representation of  $\overline \Lambda$. The details are left to the reader.
 \end{proof}

The next theorem generalizes \cite[Theorem 4.3]{Stan}.
 
 \begin{theorem} Let $S,T \subset [n]$. The enhanced crossing and nesting numbers $\eCr(\overline \Lambda)$ and $\eNe(\overline \Lambda)$ have a symmetric joint distribution over all $r$-colored enhanced partitions $\overline \Lambda$ of $[n]$ with 
 $ \min(\overline\Lambda) \setminus \max(\overline\Lambda) = S$ and $ \max(\overline\Lambda) \setminus \min(\overline\Lambda) = T$.
  \end{theorem}
  
  \begin{proof}
  Define $\phi$ as the map in Proposition \ref{enhanced-part-match-lem}.
One checks that if  $\overline\Lambda$ has  $ \min(\overline\Lambda) \setminus \max(\overline\Lambda) = S$ and $ \max(\overline\Lambda) \setminus \min(\overline\Lambda) = T$  then $M = \phi(\overline \Lambda)$ has
 $\min(M) =  \{\text{odd }i \in [2n]\} \cup \{ 2i : i \in S\}$ and $\max(M) =  \{ \text{even }i \in [2n] \} \cup \{ 2i -1: i \in T\}$,
 so we may invoke Lemma \ref{general-cor}.
  \end{proof}
  

While on this topic we mention a  result obtained in a similar way which does not seem to  be noted in the literature.
Here, we let 
$\overline{\NCN}_{j,k}(n,r) = \overline{\NCN}_{k,j}(n,r)$ 
 denote the number of $r$-colored enhanced partitions $\overline\Lambda$ of $[n]$ with $\eCr(\overline\Lambda) < j$ and $\eNe(\overline\Lambda) < k$.


\begin{proposition} For all   $j,k,n,r\in \PP$, we have $\NCN_{j,k}(n+1,r) = \sum_{i=0}^n \binom{n}{i} \overline{\NCN}_{j,k}(i,r)$.
\end{proposition}

This  statement is equivalent to 
the  identity  of exponential generating functions 
\be\label{notto}\barr{c} \frac{d}{dx} \(\sum_{n\geq 0} \NCN_{j,k}(n) \tfrac{x^n}{n!}\)= e^x\(\sum_{n\geq 0} \overline{\NCN}_{j,k}(n) \tfrac{x^n}{n!}\)\earr\ee
 which  in the case $r=1$
  provides a way of deriving either of the two main propositions in Bousquet-M\'elou and Xin's paper \cite{BousquetMelou2} from the other.
  

\def\cV{\mathcal V}
\def\cH{\mathcal H}

\begin{proof}
Consider, in slight contrast to Proposition \ref{enhanced-part-match-lem},  the map sending  an $r$-colored enhanced partition  $\overline \Lambda = (P,\varphi)$ of $[n]$ to the unique $r$-colored matching $\Lambda'=(P',\varphi') $ of $[2n+2]$ with 
$\Arc(P') = \left\{ (2i,2j+1) : (i,j) \in \overline{\Arc}(P)\right\}$ and $\varphi'(2i,2j+1) = \varphi(i,j)$. One checks that this is a bijection from  $r$-colored enhanced partitions $\overline \Lambda$ with $\eCr(\overline \Lambda) < j$ and $\eNe(\overline \Lambda) < k$ to $r$-colored $j$-noncrossing $k$-nonnesting matchings $M$ of $[2n+2]$ such that 
\ben
\item[(i)] $2i \in \min(M)$ and $2i-1 \in \max(M)$ for all $i \in [n+1]$;
\item[(ii)] $2i$ and $2i+1$ are never both isolated points of $M$ 
 for any $i \in [n]$.
\een 
Lemma \ref{part-match-lem} (which is stated in terms of uncolored partitions but extends easily to the colored case) affords a bijection from $r$-colored $j$-noncrossing $k$-nonnesting partitions of $[n+1]$ to the set of $r$-colored $j$-noncrossing $k$-nonnesting matchings $M$ of $[2n+2]$  satisfying only condition (i). Given this observation the proposition  follows by a basic counting argument. 
 \end{proof}

Before proceeding to our next topic, we observe a few corollaries of this proposition. 

\begin{corollary}
For any $j,k,r \in \PP$ the formal power series 
 $\sum_{n\geq 0} \overline{\NCN}_{j,k}(n,r) x^n $ is  rational.
\end{corollary}

\begin{proof}
This follows from Corollary \ref{rat-cor} since the proposition shows that if $F(x) = \sum_{n\geq 0} \NCN_{j,k}(n,r)x^n$ then $\sum_{n\geq 0} \overline{\NCN}_{j,k}(n,r) x^n = \frac{1}{x}\(F(\frac{x}{x+1})-1\)$.
\end{proof}

Let $\NC_k(n,r)$ and $\overline{\NC}_k(n,r)$ denote the number of $k$-noncrossing (equivalently, $k$-nonnesting) $r$-colored partitions of $[n]$. The next corollary is immediate from \eqref{notto}:

\begin{corollary} If either of the formal power series $\sum_{n\geq 0} \NC_k(n,r) x^n$ or $\sum_{n\geq 0} \overline{\NC}_k(n,r) x^n$ is D-finite then both are.
\end{corollary}

In particular, let $\overline C_n(r) = \overline {\NC}_2(n,r)$. By Theorem \ref{intro-prop} we then have:

\begin{corollary} The formal power series $\sum_{n\geq 0} \overline C_n(2) x^n$ is D-finite.
\end{corollary}

It is straightforward but not very instructive to derive an exact polynomial recurrence for $\overline C_n(r)$ from Theorem \ref{intro-prop}, and  we omit these details.

\subsection{Permutations}

Corteel \cite{Corteel} first introduced crossings and nestings for permutations (i.e., bijections $[n] \to [n]$). Burill, Mishna, and Post \cite{BMP} extended Corteel's notion to defines $k$-crossings and $k$-nestings in permutations, and Yen \cite{Yen} has recently considered such crossings and nestings in colored permutations.  Here we connect some of the results in \cite{Yen} to our methods here.

Given a permutation $\sigma$ of $[n]$:
\begin{itemize}
\item Let $A_\sigma^+$ denote the set $\{ (i, \sigma(i)) : i \in [n] \text{ such that }i\leq \sigma(i) \}$. 
\item Let $A_\sigma^-$ denote the set $\{ (\sigma(i), i) : i \in [n] \text{ such that }\sigma(i) < i\}$.
\end{itemize}
Following Yen, an \emph{$r$-colored permutation} of $[n]$ is a triple $\Sigma = (\sigma, \varphi^+, \varphi^-)$ where $\sigma$ is a permutation of $[n]$ and $\varphi^{\pm}$ are maps  $A_\sigma^{\pm} \to [r]$. The standard representation of such a triple is given by drawing $n$ dots $\bullet$ in increasing order on a horizontal line, with the  arcs $A_\sigma^+$ and $A_\sigma^-$ labeled by $\varphi^\pm$ drawn in the upper and lower half plane respectively to connect the corresponding dots (counted as 1,2,\dots,$n$ from left to right).  For example, if $a,b,c,d,e,f,g,h,i \in [r]$ then
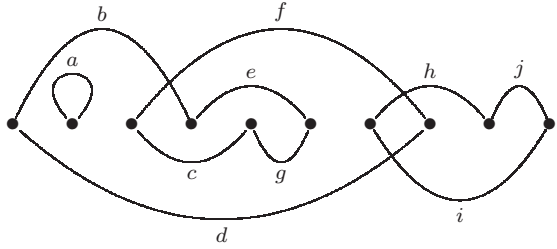
\[ 
 \xy<0.0cm,0.0cm> \xymatrix@R=-0.0cm@C=0.6cm{
*{\bullet} 
\ar @/^3.0pc/ @{-} [rrr]^b  
\ar @/_3.0pc/ @{-} [rrrrrrr]_d & 
*{\bullet}  
\ar @(ul,ur) @{-}^a &
*{\bullet}    
\ar @/^3.0pc/ @{-} [rrrrr]^f
\ar @/_1.2pc/ @{-} [rr]_c
  &
*{\bullet}   
\ar @/^1.2pc/ @{-} [rr]^e &
*{\bullet} 
\ar @/_1.2pc/ @{-} [r]_g&
*{\bullet} &
*{\bullet}    
\ar @/^1.2pc/ @{-} [rr]^h
\ar @/_2.4pc/ @{-} [rrr]_i  &
*{\bullet}  &
*{\bullet}  
\ar @/^1.2pc/ @{-} [r]^j &
*{\bullet}  
}\endxy
\]
is the standard representation of an $r$-colored permutation $\Sigma$ with $\sigma = (1,4,6,5,3,8)(2)(7,9,10)$ (written in cycle notation). Observe that $(A_\sigma^-,\varphi^-)$ represents the same data as an $r$-colored partition of $[n]$, while $(A_\sigma^+,\varphi^+)$ is an object intermediate between a colored partition and a colored enhanced partition.
We mention that the number of $r$-colored permutations of $[n]$ is $r^n n!$, and that these objects are naturally identified with elements of the wreath product $ (\ZZ/r \ZZ) \wr S_n$.

An $r$-colored permutation $\Sigma = (\sigma, \varphi^+, \varphi^-)$ has a $k$-crossing (respectively, $k$-nesting) if either of the following holds:
\begin{itemize}
\item There is a sequence of arcs $\{ (i_t,j_t) \}_{t \in [k]}\subset A_\sigma^+$, all labeled by $\varphi^+$ with the same color, which form an enhanced $k$-crossing (respectively, enhanced $k$-nesting) in the sense of \eqref{enhanced-cross-nest}.
\item There is a sequence of arcs $\{ (i_t,j_t) \}_{t \in [k]}\subset A_\sigma^-$, all labeled by $\varphi^-$ with the same color, which form a $k$-crossing (respectively, $k$-nesting) in the usual sense of \eqref{cross-nest}.
\end{itemize}
Let $\Cr(\Sigma)$ and $\Ne(\Sigma)$ be the largest integers $k$ such that $P$ has a $k$-crossing or $k$-nesting, respectively. Say that $\Sigma$ is $j$-noncrossing or $k$-nonnesting if $\Cr(\Sigma)<j$ or $\Ne(\Sigma) < k$.

The following result does not appear to be noted in the literature, and allows us to give an alternate proof of theorems in \cite{BMP} and \cite{Yen} directly from Lemma \ref{general-cor}.

\begin{proposition}\label{permutation-lem}
There is a bijection from  $r$-colored $j$-noncrossing $k$-nonnesting permutations of $[n]$ to  pairs $( \Lambda^+,\Lambda^-)$ of  $r$-colored $j$-noncrossing $k$-nonnesting  matchings of $[2n]$
such that
\begin{itemize}
\item[(a)] $\{1,3,5,\dots,2n-1\}$ is the disjoint union of $ \min(\Lambda^+) \setminus \max(\Lambda^+)$ and $\max(\Lambda^-)\setminus \min(\Lambda^-)$;
\item[(b)] $\{2,4,6,\dots,2n\}$ is the disjoint union of $\max(\Lambda^+) \setminus \min(\Lambda^+)$ and $ \min(\Lambda^-) \setminus\max(\Lambda^-)$.
\end{itemize}
\end{proposition}

\begin{proof}
Given an $r$-colored permutation $\Sigma = (\sigma,\varphi^+,\varphi^-)$ of $[n]$, let  $\Lambda^\pm = (P^\pm, \phi^\pm)$ be the unique $r$-colored partitions of $[2n]$ such that 
\begin{itemize}
\item $\Arc(P^+) = \{(2i-1,2j) : (i,j) \in A_\sigma^+\}$ and  $\phi^+(2i-1,2j) = \varphi^+(i,j)$ for $(i,j) \in A_\sigma^+$;
\item $\Arc(P^-) = \{(2i,2j-1) : (i,j) \in A_\sigma^-\}$ and  $\phi^-(2i,2j-1) = \varphi^-(i,j)$ for $(i,j) \in A_\sigma^-$.
\end{itemize}
Noting Lemma \ref{part-match-lem} and Proposition \ref{enhanced-part-match-lem}, it is straightforward to check that $\Sigma$ is $j$-noncrossing and $k$-nonnesting if and only if both $\Lambda^\pm$ are, and that $\Lambda^\pm$ are colored matchings satisfying conditions (a) and (b).

To construct an inverse to the map $\Sigma \mapsto (\Lambda^+,\Lambda^-)$, 
fix a pair $(\Lambda^+,\Lambda^-)$ of $r$-colored matchings  of $[2n]$ satisfying (a) and (b). Define  $\sigma : [n] \to [n]$ by setting $\sigma(i)$ equal to the unique  $j \in [n]$ such that  $(2i-1,2j)$ is an arc of $\Lambda^+$ or  $(2j,2i-1)$ is an arc of $\Lambda^-$;  conditions (a) and (b) ensure that exactly one of these cases occurs, and hence also that the map $\sigma$ is a permutation. Form an $r$-colored permutation $\Sigma = (\sigma,\varphi^+,\varphi^-)$ by coloring $(i,\sigma(i))$ or $(\sigma(i),i)$ in the same way as the corresponding arc $(2i-1,2j)$ or $(2j,2i-1)$ in $\Lambda^\pm$. By construction the map $(\Lambda^+,\Lambda^-) \mapsto \Sigma$  given here and the map  $\Sigma \mapsto (\Lambda^+,\Lambda^-)$  given in the previous paragraph are inverses of each other.
\end{proof}

As one application of the preceding proposition, we state the following theorem. Yen proves a more detailed version of this statement by different methods in \cite{Yen}.

\begin{theorem}[See Yen \cite{Yen}]
The  crossing and nesting numbers $\Cr(\Sigma)$ and $\Ne(\Sigma)$ have a symmetric joint distribution over all $r$-colored permutations of $[n]$. 
\end{theorem}

\begin{proof}

Lemma \ref{general-cor} implies that the numbers $\Cr(\Lambda^+,\Lambda^-) \omdef= \max\{ \Cr(\Lambda^\pm)\}$
and 
$\Ne(\Lambda^+,\Lambda^-) \omdef= \max\{ \Ne(\Lambda^\pm)\}$
have a symmetric joint distribution over all 
 pairs $(\Lambda^+,\Lambda^-)$ of $r$-colored matchings of $[2n]$ for which the four sets $\min(\Lambda^\pm) \setminus \max(\Lambda^\pm) $ and $\max(\Lambda^\pm) \setminus \min(\Lambda^\pm) $ are fixed.
  The theorem therefore follows by Proposition \ref{permutation-lem}.
\end{proof}

Let $\NCN^{Permute}_{j,k}(n,r)$ (respectively, $\NC^{Permute}_k(n,r)$) denote the number of $j$-noncrossing $k$-nonnesting (respectively, $k$-noncrossing) $r$-colored permutations of $[n]$.
Yen proves the following in \cite[Section 4]{Yen}:
\begin{corollary}[See Yen \cite{Yen}]
For all $j,k,r \in \PP$ the  power series $\sum_{n\geq 0} \NCN^{Permute}_{j,k}(n,r) x^n$ is rational.
\end{corollary}

The formal power series $\sum_{n\geq 0} \NC^{Permute}_k(n,r) x^n$ is algebraic when $k=2$ and $r=1$ since, as noted in \cite[Table 2]{BMP}, the number 
of noncrossing uncolored permutations of $[n]$ is again given by the $n$th Catalan number; i.e.,
$\NC^{Permute}_2(n,1) =  \frac{1}{n+1} \binom{2n}{n}.$
However, we are left with this question, still open even in the case $k=r=2$.

\begin{question}
For which $k,r \in \PP$ is the power series $\sum_{n\geq 0} \NC^{Permute}_k(n,r) x^n$ D-finite? 
\end{question}

Finally, we mention that there are also notions of crossings and nestings which been studied for permutations of type B; see \cite{Hamdi}.
Type B permutations are in bijection with 2-colored permutations, and it would be interesting to know of any connection between the results in \cite{Hamdi} and \cite{Yen}.

\subsection{Tangled diagrams}

Chen, Qin, and Reidys introduced tangled diagrams in \cite{Tangle} as a combinatorial framework for efficient prediction algorithms involved in  interactions between RNA molecules. Tangled diagrams a further studied, for example, if \cite{Tangle2}. Slightly generalizing the construction in \cite{Tangle}, we define an \emph{$r$-colored tangled diagram} on $[n]$ as a labeled graph on the vertices $1,2,\dots,n$ (drawn in increasing order on a horizontal line) with arcs labeled by $[r]$  (drawn in the upper half plane between vertices), such that at any vertex at most two arcs  meet in one of the following local configurations:
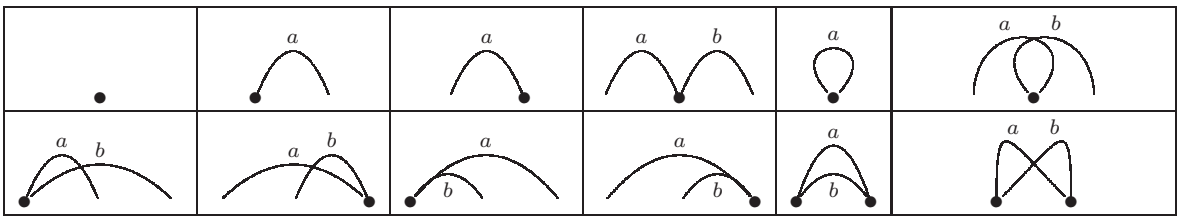
\[
\barr{|c|c|c|c|c|c|} \hline
 \xy<0.0cm,0.0cm> \xymatrix@R=-0.0cm@C=0.8cm{
*{\bullet}  
}\endxy
&
 \xy<0.0cm,0.0cm> \xymatrix@R=-0.0cm@C=0.8cm{
*{\bullet}  \ar @/^1.5pc/ @{-} [r]^a
&   
}\endxy
&
 \xy<0.0cm,0.0cm> \xymatrix@R=-0.0cm@C=0.8cm{
& 
*{\bullet}  \ar @/_1.5pc/ @{-} [l]_a 
}\endxy
&
\xy<0.0cm,0.0cm> \xymatrix@R=-0.0cm@C=0.8cm{
 &
*{\bullet}  \ar @/_1.5pc/ @{-} [l]_a \ar @/^1.5pc/ @{-} [r]^b
   & 
}\endxy
&
 \xy<0.0cm,0.0cm> \xymatrix@R=-0.0cm@C=0.8cm{
*{\bullet}  \ar @(ul,ur) @{-}^a
}\endxy
&
 \xy<0.0cm,0.0cm> \xymatrix@R=0.2cm@C=0.6cm{
&&
*{\bullet}  
 \ar @(ul,u) @{-} [r]^{\hspace{6mm}b}
  \ar @(ur,u) @{-} [l]_{\hspace{-6mm}a}
   & &
}\endxy
\\ \hline
 \xy<0.0cm,0.0cm> \xymatrix@R=-0.0cm@C=0.8cm{
*{\bullet}  \ar @/^1.2pc/ @{-} [rr]^b \ar @/^1.5pc/ @{-} [r]^a
   & &
}\endxy
&
 \xy<0.0cm,0.0cm> \xymatrix@R=-0.0cm@C=0.8cm{
 & & 
*{\bullet}  \ar @/_1.2pc/ @{-} [ll]_a \ar @/_1.5pc/ @{-} [l]_b
}\endxy
&
 \xy<0.0cm,0.0cm> \xymatrix@R=-0.0cm@C=0.8cm{
*{\bullet}  \ar @/^1.5pc/ @{-} [rr]^a \ar @/^0.9pc/ @{-} [r]_b
   & &
}\endxy
&
 \xy<0.0cm,0.0cm> \xymatrix@R=-0.0cm@C=0.8cm{
 & & 
*{\bullet}  \ar @/_1.5pc/ @{-} [ll]_a \ar @/_0.9pc/ @{-} [l]^b
}\endxy
&
 \xy<0.0cm,0.0cm> \xymatrix@R=-0.0cm@C=0.8cm{
*{\bullet} & 
*{\bullet}  \ar @/_1.8pc/ @{-} [l]_a   \ar @/_0.9pc/ @{-} [l]^b   
}\endxy
&
 \xy<0.0cm,0.0cm> \xymatrix@R=-0.0cm@C=0.8cm{
*{\bullet} \ar @(ur,u) @{-} [r]^b & 
*{\bullet}  \ar @(ul,u) @{-} [l]_a   
}\endxy
\\
\hline
\earr
\]
Here the labels $a,b \in [r]$ are arbitrary. While the vertices in these diagrams are labeled by $1,2,\dots,n$ from left to right,  we usually omit these labels from our drawings.
For example, 
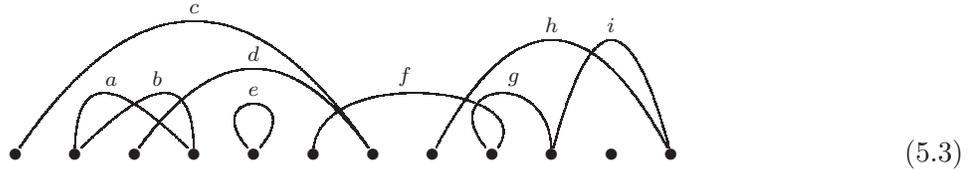
\be\label{tangled-ex} 
 \xy<0.0cm,0.0cm> \xymatrix@R=0.2cm@C=0.6cm{
*{\bullet}
 \ar @/^4.2pc/ @{-} [rrrrrr]^c&
*{\bullet} 
\ar @(ur,u) @{-} [rr]^b &
*{\bullet} 
\ar @/^2.7pc/ @{-} [rrrr]^d& 
*{\bullet}  
\ar @(ul,u) @{-} [ll]_a &
*{\bullet}  
\ar @(ul,ur) @{-}^e&
*{\bullet}&
*{\bullet} &
*{\bullet}\ar @/^3.6pc/ @{-} [rrrr]^h&
*{\bullet}  
 \ar @(ul,u) @{-} [r]^{\hspace{6mm}g}
  \ar @(ur,u) @{-} [lll]_{\hspace{-6mm}f}& 
*{\bullet}
\ar @/^3.6pc/ @{-} [rr]^i&
*{\bullet}&
*{\bullet}
}\endxy
\ee
is  an $r$-colored tangled diagram on $[n]$ with $n=12$ (here the letters $a,b,c,\dots$ indicate arbitrary elements of $[r]$).
The notion of tangled diagrams in \cite{Tangle} coincides precisely with 1-colored tangled diagrams in the sense just given.

To each $r$-colored tangled diagram $T$ on $[n]$ one associates an $r$-colored matching $\eta(T)$ of $[2n]$, called the \emph{inflation} of $T$, by doubling  each vertex  such that the local configurations of arcs shown above are respectively transformed  to the configurations below:
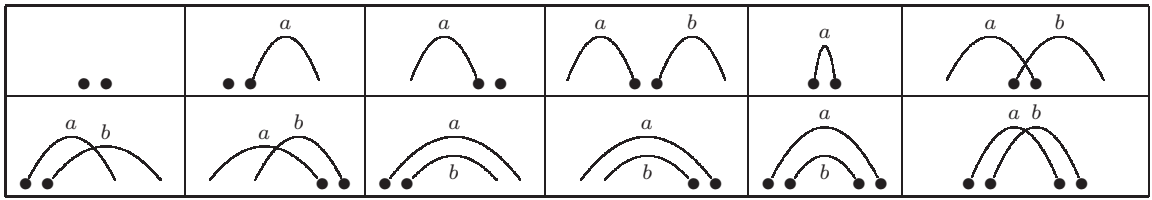
\[
\barr{|c|c|c|c|c|c|} \hline
 \xy<0.0cm,0.0cm> \xymatrix@R=-0.0cm@C=0.1cm{
*{\bullet} & *{\bullet} 
}\endxy
&
 \xy<0.0cm,0.0cm> \xymatrix@R=-0.0cm@C=0.1cm{
*{\bullet}  & *{\bullet}  \ar @/^1.5pc/ @{-} [rrr]^a
&   & &
}\endxy
&
 \xy<0.0cm,0.0cm> \xymatrix@R=-0.0cm@C=0.1cm{
& &&
*{\bullet}  \ar @/_1.5pc/ @{-} [lll]_a  & *{\bullet}  
}\endxy
&
\xy<0.0cm,0.0cm> \xymatrix@R=-0.0cm@C=0.1cm{
 & &&
*{\bullet}  \ar @/_1.5pc/ @{-} [lll]_a  & *{\bullet}  \ar @/^1.5pc/ @{-} [rrr]^b
   & & &
}\endxy
&
 \xy<0.0cm,0.0cm> \xymatrix@R=-0.0cm@C=0.1cm{
*{\bullet}  \ar @/^1.2pc/ @{-} [r]^a & *{\bullet} 
}\endxy
&
\xy<0.0cm,0.0cm> \xymatrix@R=-0.0cm@C=0.1cm{
 & &&
*{\bullet}  \ar @/^1.5pc/ @{-} [rrrr]^b  & *{\bullet} \ar @/_1.5pc/ @{-} [llll]_a 
   & &&
}\endxy
\\ \hline
 \xy<0.0cm,0.0cm> \xymatrix@R=-0.0cm@C=0.1cm{
*{\bullet} \ar @/^1.5pc/ @{-} [rrrr]^a & *{\bullet}  \ar @/^1.2pc/ @{-} [rrrrr]^b 
   & &&&&
}\endxy
&
 \xy<0.0cm,0.0cm> \xymatrix@R=-0.0cm@C=0.1cm{
&   & &&& *{\bullet}  \ar @/_1.2pc/ @{-} [lllll]_a  & *{\bullet} \ar @/_1.5pc/ @{-} [llll]_b 
}\endxy
&
 \xy<0.0cm,0.0cm> \xymatrix@R=-0.0cm@C=0.1cm{
*{\bullet} \ar @/^1.5pc/ @{-} [rrrrrr]^a & *{\bullet}  \ar @/^0.9pc/ @{-} [rrrr]_b 
   & &&&&
}\endxy
&
 \xy<0.0cm,0.0cm> \xymatrix@R=-0.0cm@C=0.1cm{
&   & &&& *{\bullet}  \ar @/_0.9pc/ @{-} [llll]^b &  *{\bullet} \ar @/_1.5pc/ @{-} [llllll]_a
}\endxy
&
 \xy<0.0cm,0.0cm> \xymatrix@R=-0.0cm@C=0.1cm{
*{\bullet} &
*{\bullet} & & &
*{\bullet}    \ar @/_0.9pc/ @{-} [lll]^b   &
*{\bullet}  \ar @/_1.8pc/ @{-} [lllll]_a  
}\endxy
&
 \xy<0.0cm,0.0cm> \xymatrix@R=-0.0cm@C=0.1cm{
&& *{\bullet} &
*{\bullet} & &&
*{\bullet}    \ar @/_1.8pc/ @{-} [llll]_a   &
*{\bullet}  \ar @/_1.8pc/ @{-} [llll]_b &&
}\endxy
\\
\hline
\earr
\]
For example, is $T$ is the diagram in \eqref{tangled-ex} then $\eta(T)$ is given by
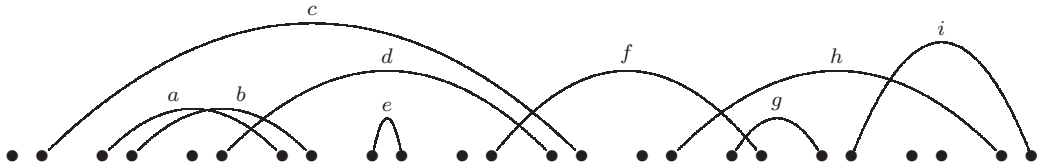
\[
 \xy<0.0cm,0.0cm> \xymatrix@R=0.2cm@C=0.2cm{
*{\bullet}
&
*{\bullet} \ar @/^4.2pc/ @{-} [rrrrrrrrrrrrrrrrrr]^c&&
*{\bullet} 
 \ar @/^1.5pc/ @{-} [rrrrrr]^{\hspace{-5mm}a} &
*{\bullet}
 \ar @/^1.5pc/ @{-} [rrrrrr]^{\hspace{5mm}b}&&
*{\bullet} & 
*{\bullet}\ar @/^2.7pc/ @{-} [rrrrrrrrrrr]^d&&
*{\bullet}&
*{\bullet}&&
*{\bullet}  
\ar @/^1.2pc/ @{-} [r]^e&
*{\bullet}&&
*{\bullet}&
*{\bullet}
\ar @/^2.7pc/ @{-} [rrrrrrrrr]^f&&
*{\bullet} 
&
*{\bullet}&&
*{\bullet}&
*{\bullet}
\ar @/^2.7pc/ @{-} [rrrrrrrrrrr]^h&&
*{\bullet}
\ar @/^1.2pc/ @{-} [rrr]^g&
*{\bullet}  
& &
*{\bullet}&
*{\bullet}\ar @/^3.6pc/ @{-} [rrrrrr]^i&&
*{\bullet}&
*{\bullet}&&
*{\bullet}&
*{\bullet}
}\endxy
 \]
Define the crossing and nesting numbers of $\Cr(T)$ and $\Ne(T)$ of a tangled diagram $T$ to be the numbers $\Cr(\eta(T))$ and $\Ne(\eta(T))$ respectively.
 In the preceding example, we have $\Cr(T) = \Ne(T) = 1$ if the arc colors are all distinct, and $\Cr(T) = \Ne(T) = 3$ if the arc colors all coincide. 
 
 \begin{remark} Our definition of the inflation map $\eta$ differs slightly when $r=1$ from the one given in \cite[Section 2.2]{Tangle}, but the difference is only that inflation here inserts more isolated points into the tangled diagram $T$ to form $\eta(T)$. Hence,  the crossing and nesting numbers in \cite{Tangle}, which  are given in exactly  the same way as the crossing and nesting numbers of the corresponding inflated diagram, coincide with our definition of $\Cr(T)$ and $\Ne(T)$.
 \end{remark}
 
The following proposition highlights the utility of the inflation map $\eta$ as we have defined it.

\begin{proposition}\label{tangled-lem} The inflation map $\eta$ defines a bijection from the set of $r$-tangled diagrams on $[n]$ to the subset of $r$-colored matchings $M$ on $[2n]$ such that for each $i \in [n]$ if  $2i$ is an isolated point of $M$ then $2i-1 \in \max(M)$ and if $2i-1$ is an isolated point of $M$ then $2i \in \min(M)$.
\end{proposition}

\begin{proof}[Proof Sketch] 
The proposition is intuitively clear since tangled diagrams and matchings are completely determined by the local arc configurations in their standard representations, and since 
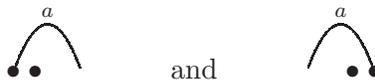
\[ \xy<0.0cm,0.0cm> \xymatrix@R=-0.0cm@C=0.1cm{
*{\bullet} \ar @/^1.5pc/ @{-} [rrr]^a & *{\bullet}  
&   & &
}\endxy
\qquand
 \xy<0.0cm,0.0cm> \xymatrix@R=-0.0cm@C=0.1cm{
& &&
*{\bullet}   & *{\bullet}   \ar @/_1.5pc/ @{-} [lll]_a
}\endxy
\]
are the only arc configurations involving  consecutive vertices $2i-1$ $2i$ (for $i \in [n]$) which can occur in a matching of $[2n]$ but which cannot occur in the inflation $\eta(T)$ of a tangled diagram on $[n]$. 
\end{proof}

Taking $\phi = \eta$ in Lemma \ref{general-cor} thus shows:

\begin{theorem} The  crossing and nesting numbers $\Cr(T)$ and $\Ne(T)$ have a symmetric joint distribution over all $r$-colored tangled diagrams $T$ on $[n]$.
\end{theorem} 

In turn, combining Proposition \ref{tangled-lem} with Theorem \ref{summary-thm} gives us this generalization of \cite[Theorems 3.6 and 3.7]{Tangle}.

\begin{theorem} There is a bijection
$T \mapsto (\bflambda^0,\bflambda^1,\dots,\bflambda^{2n})$ from the set of $r$-colored  tangled diagrams  on $[n]$ to the subset of semi-oscillating $r$-partite tableaux of length $2n$ such that
for each $i \in [n]$ the transition from $\bflambda^{2i-2}$ to $\bflambda^{2i-1}$ to  $\bflambda^{2i}$ does not   consist of either 
\begin{itemize}
\item Adding a box then doing nothing;
\item Doing nothing then deleting a box.
\end{itemize}
Further, $T$ is $j$-noncrossing and $k$-nonnesting if and only if  the $r$ components of each $\bflambda^i$ all have fewer than $j$ rows and and fewer than $k$ columns.
\end{theorem}

\begin{proof}
The bijection from $r$-colored matchings of $[2n]$ to semi-oscillating $r$-partite tableaux of length $2n$ described in Theorem \ref{summary-thm} is achieved explicitly by first applying the map in Lemma \ref{intro-lem} to an $r$-colored matching and then applying the map in Theorem \ref{match-map} to resulting $r$ uncolored components. With respect to the resulting bijection $M \mapsto (\bflambda^0,\bflambda^1,\dots,\bflambda^n)$, one checks that $i \in \min(M)$ if and only if $\bflambda^{i-1} \subseteq \bflambda^i$ and $i \in \max(M)$ if and only if $\bflambda^{i-1} \supseteq \bflambda^i$. The theorem here follows by combining this observation with Proposition \ref{tangled-lem}. 
\end{proof}

Let $\NCN^{Tangle}_{j,k}(n,r)$ (respectively, $\NC^{Tangle}_k(n,r)$) denote the number of $j$-noncrossing $k$-nonnesting (respectively, $k$-noncrossing) $r$-colored tangled diagrams on $[n]$.
Given the preceding result, the following corollary (like Corollary \ref{matching-rat}) follows by arguments  similar to the proofs of Theorem \ref{graph-thm} and Corollary \ref{rat-cor}.

\begin{corollary}
For all $j,k,r \in \PP$ the formal power series $\sum_{n\geq 0} \NCN^{Tangle}_{j,k}(n,r) x^n$ is rational.
\end{corollary}

The authors of \cite{Tangle2} prove that $\sum_{n\geq 0} \NC^{Tangle}_k(n,1) x^n$ is D-finite and derive an explicit asymptotic formula for the number of $k$-noncrossing (uncolored) tangled diagrams on $[n]$. It is expected that the methods of \cite{Tangle2} extend to colored tangled diagrams without great difficulty, but we relegate the working out of these details to this open question:

\begin{question}
For which $k,r \in \PP$ is the  power series $\sum_{n\geq 0} \NC^{Tangle}_k(n,r) x^n$ D-finite?
\end{question}

\end{document}